\documentclass[12pt]{article}

\usepackage[T1]{fontenc}
\usepackage{lmodern}
\usepackage{microtype}
\usepackage{indentfirst}
\usepackage{amsmath,amssymb,amsthm,mathtools}
\usepackage{booktabs}
\usepackage{float}
\usepackage{enumitem}
\usepackage{tikz}
\usetikzlibrary{arrows.meta,positioning,calc,shapes.geometric}
\usepackage{xcolor}
\usepackage[margin=1in]{geometry}
\usepackage{titlesec}
\usepackage[pagebackref]{hyperref}
\hypersetup{colorlinks=true}
\usepackage[nameinlink,noabbrev]{cleveref}
\numberwithin{equation}{section}
\usepackage{caption}
\usepackage{subcaption}
\usepackage{ulem}
\definecolor{wallred}{RGB}{150,66,66}
\definecolor{holegray}{RGB}{105,110,115}
\definecolor{holefill}{RGB}{247,248,249}

\newtheorem{theorem}{Theorem}[section]

\newtheorem{lemma}[theorem]{Lemma}
\theoremstyle{definition}
\newtheorem{definition}[theorem]{Definition}
\newtheorem{example}[theorem]{Example}

\newtheorem{remark}[theorem]{Remark}
\crefname{definition}{Definition}{Definitions}
\crefname{example}{Example}{Examples}
\crefname{theorem}{Theorem}{Theorems}
\crefname{corollary}{Corollary}{Corollaries}
\crefname{proposition}{Proposition}{Propositions}
\crefname{lemma}{Lemma}{Lemmas}
\crefname{appendix}{Appendix}{Appendix}
\crefname{figure}{Figure}{Figure}
\crefname{section}{Section}{Section}

\newcommand{\D}{\mathcal D}
\newcommand{\Dhat}{\widehat{\mathcal D}}

\newcommand{\eps}{\varepsilon}

\newcommand{\Root}{\operatorname{Root}}
\newcommand{\Leaf}{\operatorname{Leaf}}

\newcommand{\TC}{\mathcal{TC}}

\newcommand{\pref}{\mathsf{pref}}

\def\blue{\textcolor{blue}}

\def\magenta{\textcolor{magenta}}

\title{A symmetry on  Young tableaux with
walls and holes}

\author{Zhicong Lin$^{\color{blue} \dag}$, Feihu Liu$^{\color{blue} \ddag}$, and Jing Liu$^{\color{blue}\S}$\footnote{Corresponding author}
\\[2mm]
{\small $^{\color{blue} \dag, \S}$ Research Center for Mathematics and Interdisciplinary Sciences,}\\[-0.8ex]
{\small Shandong University, Qingdao 266237, P.R.~China}\\
{\small $^{\color{blue} \ddag}$ Center for Combinatorics, LPMC,}\\[-0.8ex]
{\small Nankai University, Tianjin 300071, P.R.~China}\\
{\small $^{\color{blue}  \S}$ Department of Mathematics,}\\[-0.8ex]
{\small Sungkyunkwan University, Suwon  16419, South Korea}\\
{\small {\color{blue} $^\dag$} Email address: linz@sdu.edu.cn}\\
{\small {\color{blue} $^\ddag$} Email address: liufeihu7476@163.com}\\
{\small {\color{blue} $^\S$} Email address: lsweet@mail.sdu.edu.cn}\\
}

\date{\today}

\begin{document}

\maketitle

\begin{abstract}
The Pons--Batle conjecture on tree-child networks was regarded as a specific and challenging problem in the phylogenetics community. We found an unexpected weighted symmetry for a class of  three-row  Young tableaux with walls and holes. This weighted symmetry implies the Fuchs--Yu identity, which together with the bijection of  Chang, Fuchs,  Liu,  Wallner, and Yu leads to a simple proof of  the Pons--Batle conjecture.  We construct an involutive  proof of this symmetry via introducing decorated binary networks as an intermediate structure.   We also provide an elementary proof of this weighted symmetry using recurrence relations.
\end{abstract}

\noindent
\begin{small}
{\it 2020 Mathematics subject classification}: 05A15; 05A19; 92D15.
\end{small}

\noindent
\begin{small}
{\it Keywords}: Young tableaux; Tree-child networks; Symmetry; Bijections.
\end{small}

%%%%%%%%%%%%%%%%%%%%%%%%%%%%%%%%%%%%%
\section{Introduction}%%%%%%%%%%%%%%%%%%%%%%%%%%%
\label{sec:introduction}
%%%%%%%%%%%%%%%%%%%%%%%%%%%%%%%%%%%%
A \blue{\it rooted binary phylogenetic network} on
$[n]:=\{1,\ldots,n\}$ is a connected directed acyclic graph without
parallel edges, with a unique root of indegree zero and outdegree one,
$n$ leaves of indegree one and outdegree zero bijectively labeled by
$[n]$, and all remaining vertices have indegree-outdegree pairs
$(1,2)$ or $(2,1)$, called \blue{\it tree vertices} and
\blue{\it reticulation vertices}, respectively. 
An important subclass is
formed by \blue{\it tree-child networks}~\cite{PonsBatle2021}, in which every nonleaf vertex has
 at least one child that is not a reticulation vertex. 
The tree-child condition
ensures that every vertex lies on a directed path to a leaf whose
intermediate vertices are tree vertices, while allowing reticulation
events.

Let $\TC_{n,k}$ be the set of all tree-child networks with
$n$ labeled leaves and $k$ reticulation vertices. Cardona and Zhang~\cite{CardonaZhang2020} obtained an
exact enumeration formula for $\TC_{n,k}$ by decomposing these networks into
tree-components and their associated component graphs. Another approach was developed by Fuchs, Yu, and
Zhang~\cite{FuchsYuZhang2021}, who established a word correspondence for
the maximally reticulated subclass and used it to determine the asymptotic growth of the total
number of tree-child networks. Motivated by their correspondence, Pons
and Batle~\cite{PonsBatle2021} introduced a larger family
$\mathcal C_{n,k}$ of constrained words, with $k$ letters appearing
three times and the other $n-k$ letters twice, and conjectured that
\begin{equation}\label{eq:intro-pb-words}
|\TC_{n,k}|=\frac{n!}{(n-k)!}c_{n-1,k},
\end{equation}
where $c_{n,k}:=|\mathcal C_{n,k}|$.
Since 
$c_{n,k}=c_{n,k-1}+(2n+k-1)c_{n-1,k}$, this relation  would yield a
 two-term recurrence for $|\TC_{n,k}|$:
\begin{equation}\label{rec:TCNs}
|\TC_{n,k}|=(n-k+1)|\TC_{n,k-1}|+\frac{n(2n+k-3)}{n-k}|\TC_{n-1,k}|.
\end{equation}
Fuchs, Liu, and Yu~\cite{FuchsLiuYu2021} proved~\eqref{eq:intro-pb-words} bijectively for 
the important class of one-component tree-child networks. 

The Pons--Batle conjecture also admits a formulation in terms of \blue{\it Young
tableaux with walls}, introduced by Banderier, Marchal, and
Wallner~\cite{BanderierMarchalWallner2018} and further studied by
Banderier and Wallner~\cite{BanderierWallner2021}. A wall indicates
that the usual increasing condition for Young tableaux
 is omitted across the corresponding
pair of adjacent cells. Let $a_{n,k}$ count three-row Young tableaux with row
lengths $(k,n,n)$, from top to bottom, all left-justified and with walls
between consecutive cells in the bottom row. Let $b_{n,k}$ count the analogous Young
tableaux having $n$ cells in each of the first two rows and $k$ bottom
cells in arbitrary columns, with walls
between consecutive cells in the bottom row. In both families, the first two rows are increasing from left to right and
each column are increasing from bottom to top, but distinct bottom entries need not be
increasing. In fact, the first class of Young
tableaux with walls of row lengths $(k,n,n)$ are in natural bijection with the constrained words $\mathcal C_{n,k}$, and thus 
$a_{n,k}=c_{n,k}$.
Fuchs and Yu conjectured (see~\cite[Theorem~1.1]{LinEtAl2026}) the following relationship between $a_{n,k}$ and $b_{n,k}$:
\begin{equation}\label{eq:intro-fuchs-yu}
2^{n-k}a_{n,k}=(n-k+1)!\,b_{n,k}.
\end{equation}
In view of $a_{n,k}=c_{n,k}$ and the bijective proof of $2^{n-k}|\TC_{n+1,k}|=(n+1)!b_{n,k}$ by Chang, Fuchs, Liu, Wallner, and Yu~\cite[Theorem~3.4]{ChangEtAl2024},  the Fuchs--Yu identity~\eqref{eq:intro-fuchs-yu} is equivalent to the Pons--Batle conjecture stated  in~\eqref{eq:intro-pb-words}.

Over the past few years,  the Pons--Batle conjecture~\eqref{eq:intro-pb-words} has been regarded as a specific and challenging problem
in the phylogenetics community.  
Until recently, Lin, Liu, Liu, Liu, and Xin~\cite{LinEtAl2026} proved~\eqref{eq:intro-fuchs-yu} using recurrences and algebraic
identities, thereby settling the Pons--Batle conjecture~\eqref{eq:intro-pb-words}. In pursuing a combinatorial understanding of~\eqref{eq:intro-pb-words} and~\eqref{eq:intro-fuchs-yu}, 
Liu, Wallner, and Yu~\cite{LiuWallnerYu2026} introduced a
 class of  three-row  Young tableaux with walls and holes. Although they were unable to give a combinatorial proof of~\eqref{eq:intro-fuchs-yu}, they did verify~\eqref{eq:intro-fuchs-yu}  for all $0\le k\le250$ and arbitrary $n$ using lattice-path generating functions and computer
algebra. More recently, Yu and
Zhang~\cite{YuZhang2026} gave a direct combinatorial proof of~\eqref{rec:TCNs} by decomposing  tree-child networks at their roots, thereby providing an alternative approach to~\eqref{eq:intro-pb-words}.

While two distinct proofs of the Pons–Batle conjecture are already known, a combinatorial proof of~\eqref{eq:intro-fuchs-yu} remains elusive.
The main objective of the present paper is to give a bijective proof of a weighted symmetry for a class of three-row Young tableaux with walls and holes, which generalizes the Fuchs--Yu identity~\eqref{eq:intro-fuchs-yu}.

\begin{definition}\label{def:wall-hole}
Let $n\ge0$ and $0\le k,\ell\le n$. Start with a diagram having $k$ cells
in the top row, $n$ cells in the middle row, and $\ell$ cells in arbitrary
positions of the bottom row, with all cells in the top and middle rows
left-justified. The unoccupied bottom positions are called \blue{\it holes},
and a wall is placed between two consecutive bottom positions unless
both are holes. A \blue{\it Young tableau with walls and holes of type $(n,k,\ell)$}
is a bijective filling of this diagram with $[n+k+\ell]$ in which
the top and middle rows increase from left to right and each column
increases from bottom to top. No increasing condition is imposed between
distinct bottom cells. Let $\D_{n,k,\ell}$ denote the set of all these tableaux and write
$d_{n,k,\ell}=|\D_{n,k,\ell}|$. 
\end{definition}

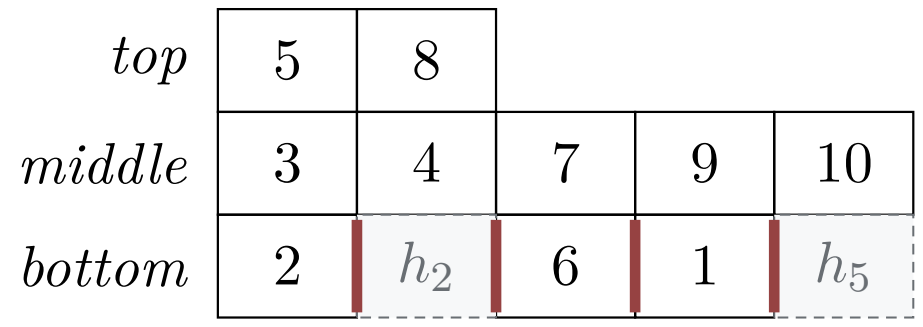
\begin{figure}
    \centering
\begin{tikzpicture}[
    x=0.92cm,
    y=0.68cm,
    line width=0.42pt,
    every node/.style={font=\small,inner sep=1pt},
    hole/.style={
        dash pattern=on 2pt off 1.4pt,
        draw=holegray,
        fill=holefill,
        line width=0.38pt
    },
    wall/.style={
        draw=wallred,
        line width=2pt
    }
]

% top row
\foreach \c/\val in {1/5,2/8}{
    \draw (\c-1,2) rectangle (\c,3);
    \node at (\c-.5,2.5) {$\val$};
}

% middle row
\foreach \c/\val in {1/3,2/4,3/7,4/9,5/10}{
    \draw (\c-1,1) rectangle (\c,2);
    \node at (\c-.5,1.5) {$\val$};
}

% occupied bottom positions
\foreach \c/\val in {1/2,3/6,4/1}{
    \draw (\c-1,0) rectangle (\c,1);
    \node at (\c-.5,.5) {$\val$};
}

% holes
\foreach \c/\lab in {2/h_2,5/h_5}{
    \draw[hole] (\c-1,0) rectangle (\c,1);
    \node[text=holegray] at (\c-.5,.5) {$\lab$};
}

% walls
\foreach \x in {1,2,3,4}{
    \draw[wall] (\x,.05) -- (\x,.95);
}

% row labels
\node[anchor=east,font=\small\itshape] at (-.18,2.5) {top};
\node[anchor=east,font=\small\itshape] at (-.18,1.5) {middle};
\node[anchor=east,font=\small\itshape] at (-.18,.5) {bottom};

\end{tikzpicture}
\caption{A Young tableau with walls and holes  of type $(5,2,3)$. \label{fig:tableau-523}}
\end{figure}

For $T\in\D_{n,k,\ell}$,
let $S(T)\subseteq[n]$ be the set of columns with occupied bottom cells.
If $j\notin S(T)$, write $h_j$ for the hole in column $j$ and put
$H(T)=\{h_j:j\notin S(T)\}$. Fig.~\ref{fig:tableau-523} shows a Young tableau with walls and holes $T\in\D_{5,2,3}$, where $S(T)=\{1,3,4\}$ and $H(T)=\{h_2,h_5\}$. 
This class of Young tableaux with walls and holes, originally introduced by Liu, Wallner, and Yu~\cite{LiuWallnerYu2026}, is motivated by the specialization
\begin{equation}\label{spec:ab}
d_{n,k,n}=a_{n,k}\quad\text{and}\quad d_{n,n,k}=b_{n,k}.
\end{equation}

Our main result is the following striking weighted symmetry for $d_{n,k,\ell}$.
\begin{theorem}\label{thm:weighted-symmetry}
For $n\ge0$ and $0\le k,\ell\le n$,
\begin{equation}\label{eq:intro-weighted-symmetry}
2^{\ell}(n-\ell+1)!\,d_{n,k,\ell}
=2^k(n-k+1)!\,d_{n,\ell,k}.
\end{equation}
\end{theorem}

Taking $\ell=n$ in Theorem~\ref{thm:weighted-symmetry} recovers the Fuchs--Yu identity~\eqref{eq:intro-fuchs-yu} in view  of~\eqref{spec:ab}. Theorem~\ref{thm:weighted-symmetry} is unexpected in the following senses: 
\begin{itemize}
\item the Fuchs--Yu identity~\eqref{eq:intro-fuchs-yu} admits such a symmetric generalization; 
\item whereas the original proof of the Fuchs--Yu identity in~\cite{LinEtAl2026} is quite involved, its symmetric generalization has a short and elementary proof;
\item this symmetric generalization possesses a natural  bijective  proof  via the introduction of decorated binary networks as an intermediate structure.
\end{itemize}
This work reflects the profound insight of  D. Hilbert: from a suitable  general standpoint, the problem becomes more accessible to investigation.

{\bf The rest of this paper is organized as follows.} In Sec.~\ref{sec:recurrence}, we provide an elementary proof of Theorem~\ref{thm:weighted-symmetry}.  Via introducing decorated binary networks as an intermediate structure, we present in Sec.~\ref{sec:bijection} a bijective proof of Theorem~\ref{thm:weighted-symmetry}.

\section{An elementary proof of Theorem~\ref{thm:weighted-symmetry}}
\label{sec:recurrence}

In this section, we prove Theorem~\ref{thm:weighted-symmetry} by induction on \(n\).
By convention, we set
\(
d_{n,k,\ell}=0
\)
unless \(n\ge 0\) and \(0\le k,\ell\le n\), and set
\(
d_{0,0,0}=1.
\)
Let $$G_{n,k,\ell}:=2^\ell(n-\ell+1)!\,d_{n,k,\ell}.$$
The main ingredients are two recurrences for \(G_{n,k,\ell}\),
with respect to \(k\) and \(\ell\), respectively. 

We first derive a recurrence for $G_{n,k,\ell}$ with respect to $k$.

\begin{lemma}\label{prop:k-recurrence}
For $n\geq1$ and $0\leq k,\ell\leq n$,
\begin{equation}\label{eq:k-recurrence}
G_{n,k,\ell}=G_{n,k-1,\ell}
+
(n-\ell+1)G_{n-1,k,\ell}
+2(n+k+\ell-1)G_{n-1,k,\ell-1}.
\end{equation}
\end{lemma}

\begin{proof}
By considering the position of the largest entry in any Young tableau in $\D_{n,k,\ell}$, it can be shown (see~\cite[Lemma~10]{LiuWallnerYu2026}) that 
 \[
d_{n,k,\ell}
=
d_{n,k-1,\ell}
+d_{n-1,k,\ell}
+(n+k+\ell-1)d_{n-1,k,\ell-1}.
\]
Multiplying both sides by \(2^\ell(n-\ell+1)!\) gives
\eqref{eq:k-recurrence}.
\end{proof}

As a consequence of \eqref{eq:k-recurrence}, we get the boundary values
\begin{equation}\label{eq:boundary-values}
G_{n,k,0}=\frac{(n+k)!(n-k+1)}{k!}=G_{n,0,k}
\end{equation}
for $0\le k\le n$. 
Indeed, specializing \eqref{eq:k-recurrence} to $\ell=0$ and using $G_{0,0,0}=1$, the first equality of~\eqref{eq:boundary-values}  follows by
induction on $n$. The second equality of~\eqref{eq:boundary-values}  follows similarly.

We next derive the hidden  recurrence with respect to $\ell$.

\begin{lemma}\label{lem:l-recurrence}
For $n\geq1$ and $0\leq k,\ell\leq n$,
\begin{equation}\label{eq:l-recurrence}
G_{n,k,\ell}=G_{n,k,\ell-1}
+
(n-k+1)G_{n-1,k,\ell}
+2(n+k+\ell-1)G_{n-1,k-1,\ell}.
\end{equation}
\end{lemma}

\begin{proof}
For $m\ge1$ and $0\le i,j\le m$, introduce  the defect
\[
\begin{aligned}
\Delta_{m,i,j}:={}
G_{m,i,j}-G_{m,i,j-1}
-(m-i+1)G_{m-1,i,j}-2(m+i+j-1)G_{m-1,i-1,j}.
\end{aligned}
\]
Thus \eqref{eq:l-recurrence} is equivalent to $\Delta_{n,k,\ell}=0$.
By \eqref{eq:boundary-values}, direct substitution gives
$\Delta_{n,0,\ell}=0$ and $\Delta_{n,k,0}=0$ for all $0\le k,\ell\le n$.

We proceed to prove $\Delta_{n,k,\ell}=0$ by induction on $n$. For $n=1$, the boundary
cases leave only $(k,\ell)=(1,1)$. By \eqref{eq:k-recurrence} and
\eqref{eq:boundary-values},
$G_{1,1,1}=G_{1,0,1}=G_{1,1,0}$, so
$\Delta_{1,1,1}=0$.
Let $n\ge2$ and assume that $\Delta_{n-1,i,j}=0$ for all
$0\le i,j\le n-1$. 
For \(1\le k,\ell\le n-1\), using \eqref{eq:k-recurrence} at levels
$n$ and $n-1$ and simplifying, we obtain
\[
\Delta_{n,k,\ell}-\Delta_{n,k-1,\ell}
=
(n-\ell)\Delta_{n-1,k,\ell}
+
2(n+k+\ell-2)\Delta_{n-1,k,\ell-1}.
\]
The right-hand side vanishes by the induction hypothesis. Hence
$\Delta_{n,k,\ell}=\Delta_{n,k-1,\ell}$, and therefore
$\Delta_{n,k,\ell}=\Delta_{n,0,\ell}=0$.

It remains to consider the cases $\ell=n$ and $k=n$.
If $\ell=n$, then by the zero convention
$\Delta_{n,k,n}=G_{n,k,n}-G_{n,k,n-1}$. Hence, for
$1\le k\le n-1$, subtracting \eqref{eq:k-recurrence} at $\ell=n-1$
from that at $\ell=n$ gives
\[
\Delta_{n,k,n}-\Delta_{n,k-1,n}
=
2(2n+k-2)\Delta_{n-1,k,n-1}=0.
\]
Thus $\Delta_{n,k,n}=\Delta_{n,0,n}=0$.

If $k=n$ and $1\le\ell\le n$, 
we have
\[
\begin{aligned}
\Delta_{n,n,\ell}-\Delta_{n,n-1,\ell}
={}&
\bigl(G_{n,n,\ell}-G_{n,n-1,\ell}\bigr)
-\bigl(G_{n,n,\ell-1}-G_{n,n-1,\ell-1}\bigr)\\
&\quad
+2(2n+\ell-2)
\bigl(G_{n-1,n-2,\ell}
-G_{n-1,n-1,\ell}\bigr).
\end{aligned}
\]
By \eqref{eq:k-recurrence} and the zero convention,
\(G_{n,n,j}=G_{n,n-1,j}\) for \(0\le j\le n\).
Applying the same recurrence at level \(n-1\) with \(k=n-1\) also gives
\(G_{n-1,n-1,j}=G_{n-1,n-2,j}\) for
\(0\le j\le n-1\). Hence
\(\Delta_{n,n,\ell}-\Delta_{n,n-1,\ell}=0\).
Since \(\Delta_{n,n-1,\ell}=0\) has already been established in the
preceding cases, it follows that \(\Delta_{n,n,\ell}=0\) for
\(0\le\ell\le n\). This completes the proof by induction.
\end{proof}

We are now in a position to prove Theorem~\ref{thm:weighted-symmetry}.

\begin{proof}[{\bf Proof of Theorem~\ref{thm:weighted-symmetry}}]
We prove by induction on \(n\). The case \(n=0\) is immediate. Assume that
\(G_{n-1,k,\ell}=G_{n-1,\ell,k}\) for all
\(0\le k,\ell\le n-1\).

Fix \(0\le k,\ell\le n\). If \(k=0\), then
\(G_{n,0,\ell}=G_{n,\ell,0}\) by
\eqref{eq:boundary-values}. Suppose that \(k\ge1\). Apply
\eqref{eq:k-recurrence} to \((k,\ell)\) and
\eqref{eq:l-recurrence} to \((\ell,k)\), together with the induction hypothesis, we obtain
\[
G_{n,k,\ell}-G_{n,\ell,k}
=
G_{n,k-1,\ell}-G_{n,\ell,k-1}.
\]
Repeating this relation for \(k,k-1,\ldots,1\) yields
\(
G_{n,k,\ell}-G_{n,\ell,k}
=
G_{n,0,\ell}-G_{n,\ell,0}.
\)
The right-hand side is zero by \eqref{eq:boundary-values}. Hence
\(G_{n,k,\ell}=G_{n,\ell,k}\), which completes the proof.
\end{proof}

%%%%%%%%%%%%%%%%%%%%%%%%%%%%%%%%%%%
\section{A bijective proof of Theorem~\ref{thm:weighted-symmetry}}
\label{sec:bijection} 
%%%%%%%%%%%%%%%%%%%%%%%%%%%%%%%%%%
This section is devoted to a bijective proof of Theorem~\ref{thm:weighted-symmetry}.

%%%%%%%%%%%%%%%%%%%%%%%%%%%%%%%%%
\subsection{Decorated Young tableaux with walls and holes} %%%
%%%%%%%%%%%%%%%%%%%%%%%%%%%%%%%%%

The factors $2^\ell$ and $(n-\ell+1)!$ in
\eqref{eq:intro-weighted-symmetry} suggest the following decoration.

\begin{definition}\label{def:decorated-tableau}
A \blue{\it decorated Young tableau with walls and holes of type $(n,k,\ell)$}
is a triple $(T;\varepsilon,\rho)$, where $T\in\D_{n,k,\ell}$ is a tableau,
$\varepsilon:S(T)\to\{L,R\}$ is a map, and
$\rho:H(T)\cup\{\star\}\to[n-\ell+1]$ is a bijection,
with $\star$ a formal symbol. Let $\Dhat_{n,k,\ell}$ denote
this family and let $\Dhat_{n}:=\bigcup\limits_{0\leq k,\ell\leq n}\Dhat_{n,k,\ell}$.
\end{definition}

\begin{example}\label{ex:decorated-tableau-523}
\captionsetup{hypcap=false}
Fig.~\ref{fig:decorated-tableau-523} shows a decorated tableau
$X=(T;\varepsilon,\rho)\in\Dhat_{5,2,3}$, where
\(
(\varepsilon(1),\varepsilon(3),\varepsilon(4))=(L,R,L)
\)
and
\(
(\rho(h_2),\rho(h_5),\rho(\star))=(3,1,2).
\)

\begin{figure}
  \centering
\begin{tikzpicture}[
    x=0.92cm,
    y=0.68cm,
    line width=0.42pt,
    every node/.style={font=\small,inner sep=1pt},
    hole/.style={
        dash pattern=on 2pt off 1.4pt,
        draw=holegray,
        fill=holefill,
        line width=0.38pt
    },
    wall/.style={
        draw=wallred,
        line width=2pt
    }
]

% top row
\foreach \c/\val in {1/5,2/8}{
    \draw (\c-1,2) rectangle (\c,3);
    \node at (\c-.5,2.5) {$\val$};
}

% middle row
\foreach \c/\val in {1/3,2/4,3/7,4/9,5/10}{
    \draw (\c-1,1) rectangle (\c,2);
    \node at (\c-.5,1.5) {$\val$};
}

% occupied bottom positions
% tableau entry centered; epsilon-label in the lower-right corner
\foreach \c/\val/\eps in {1/2/L,3/6/R,4/1/L}{
    \draw (\c-1,0) rectangle (\c,1);
    \node at (\c-.5,.5) {$\val$};
    \node[anchor=south east, font=\scriptsize]
        at (\c-.08,.08) {$\eps$};
}

% holes: display the images under rho
\foreach \c/\val in {2/3,5/1}{
    \draw[hole] (\c-1,0) rectangle (\c,1);
    \node[text=wallred] at (\c-.5,.5) {$\val$};
}

% walls
\foreach \x in {1,2,3,4}{
    \draw[wall] (\x,.05) -- (\x,.95);
}

% star and its image under rho
\node at (5.55,.5) {$\star$};
\node[text=wallred, anchor=west] at (5.78,.5) {$2$};

% row labels
\node[anchor=east,font=\small\itshape] at (-.18,2.5) {top};
\node[anchor=east,font=\small\itshape] at (-.18,1.5) {middle};
\node[anchor=east,font=\small\itshape] at (-.18,.5) {bottom};

\end{tikzpicture}

\caption{
A decorated Young tableau with walls and holes  of type $(5,2,3)$.\label{fig:decorated-tableau-523}}

\end{figure}
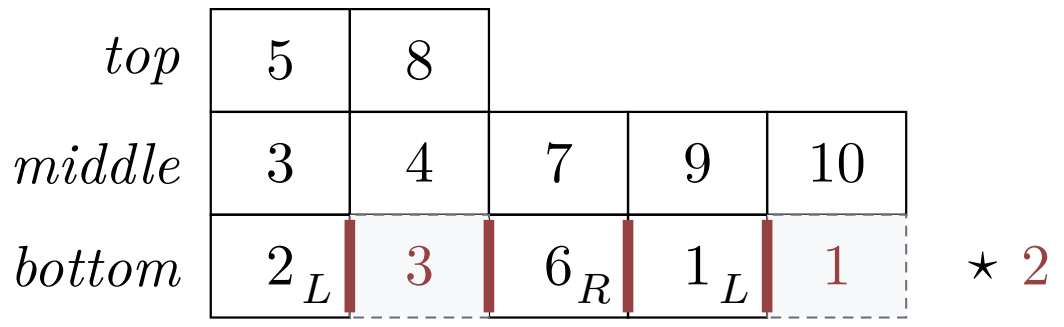
\end{example}

Since $|S(T)|=\ell$ and $|H(T)|=n-\ell$, we have
\[
|\Dhat_{n,k,\ell}|=2^{\ell}(n-\ell+1)!\,d_{n,k,\ell}.
\]
The main result of this section is the following involution on $\Dhat_n$, which proves Theorem~\ref{thm:weighted-symmetry} bijectively. 

\begin{theorem}\label{thm:bijective-symmetry}
There exists  a combinatorial involution $\Phi$ on $\Dhat_n$ such that 
$\Phi(\Dhat_{n,k,\ell})=\Dhat_{n,\ell,k}$ for all $0\leq k,\ell\leq n$. 
\end{theorem}

Our involution $\Phi$ is constructed via a correspondence with decorated binary networks that we introduce next and an
edge-reversal involution on such networks.

\subsection{Decorated binary networks and canonical postorder}\label{subsec:networks-postorder}

Let \(n\ge 0\) and \(0\le k,\ell\le n\).
A \blue{\it binary network of type \((n,k,\ell)\)} is a  directed acyclic multigraph (not necessarily connected) whose vertices consist of
\(n-k+1\) roots, \(n-\ell+1\) leaves, \(k\) tree vertices, and \(\ell\) reticulation vertices, with respective indegree-outdegree pairs
\(
(0,1), (1,0), (1,2), (2,1).
\)
Let $N$ be a binary network. 
For a vertex \(v\) of $N$, let \(E^+(v)\) and \(E^-(v)\) denote the sets of outgoing and incoming edges of \(v\), respectively. Let $\Root(N)$ and $\Leaf(N)$ be the sets of roots and leaves of $N$, respectively. Since every edge is directed from a higher vertex to a lower vertex, we suppress the directions in all our drawings. 
Note that any binary network of type \((n,k,\ell)\) has $n+k+\ell+1$ edges.

\begin{definition}\label{def:decorated-network}
A \blue{\it decorated binary network of type $(n,k,\ell)$} is a binary network $N$ of type $(n,k,\ell)$ together with the following labeling: 
\begin{enumerate}[label=\textup{(\roman*)},leftmargin=2.5em,itemsep=2pt,topsep=3pt]
\item its roots and leaves are labeled by bijections
$\rho_{\mathrm r}:\operatorname{Root}(N)\to[n-k+1]$ and
$\rho_{\mathrm \ell}:\operatorname{Leaf}(N)\to[n-\ell+1]$; 
\item for each tree vertex \(v\), its outgoing edges are labeled by a bijection
\(
\sigma_v:E^+(v)\to\{L,R\},
\)
and for each reticulation vertex \(w\), its incoming edges are labeled by a bijection
\(
\tau_w:E^-(w)\to\{L,R\}.
\)
 \end{enumerate}
\end{definition}

For the sake of convenience, we set
\(
r_i=\rho_{\mathrm r}^{-1}(i)
\)
for $i\in[n-k+1]$. 
For \(e\in E^+(v)\) (resp.,~\(e\in E^-(w)\)), we call
\(\sigma_v(e)\) (resp.,~\(\tau_w(e)\)) the \blue{\it local label of \(e\) at \(v\) (resp.,~\(w\))}.

Two decorated binary networks are \blue{\it isomorphic} if there is a directed multigraph isomorphism between their underlying binary networks  
preserving \(\rho_{\mathrm r}\), \(\rho_{\mathrm \ell}\), \(\sigma_v\), and \(\tau_w\). We denote by \(\mathcal N_{n,k,\ell}\) the set of isomorphism classes of decorated binary networks of type \((n,k,\ell)\). See Appendix~\ref{app:networks-211} for all
isomorphism classes in \(\mathcal N_{2,1,1}\).

\begin{remark}
Suppose that a tree vertex \(v\) and a reticulation vertex \(w\) are joined
by two parallel edges. Since these two edges form both \(E^+(v)\) and
\(E^-(w)\), the two local labelings determine a bijection
\[
\tau_w\circ\sigma_v^{-1}:\{L,R\}\longrightarrow\{L,R\}.
\]
Hence there are exactly two possible relative matchings,
\(\{(L,L),(R,R)\}\) and \(\{(L,R),(R,L)\}\), where each pair records the
local labels of the same edge at \(v\) and \(w\), respectively. See Appendix~\ref{app:networks-211}: (5) and (6) for an example. 
\end{remark}

The following bijection is the main ingredient in our construction of $\Phi$.

\begin{theorem}
\label{thm:tableau-network-correspondence}
For \(n\geq 0\) and \(0\leq k,\ell\leq n\), there exists  a combinatorial  bijection
\[
\Psi:
\widehat{\mathcal D}_{n,k,\ell}
\longrightarrow
\mathcal N_{n,k,\ell}.
\]
\end{theorem}

Using $\Psi$ we can now construct the involution $\Phi$ for Theorem~\ref{thm:bijective-symmetry}.

\begin{proof}[{\bf Proof of Theorem~\ref{thm:bijective-symmetry}}]
For \(N\in\mathcal N_{n,k,\ell}\), let \(\operatorname{Rev}(N)\) be
the decorated binary network obtained by reversing every edge of
\(N\). The roles of roots and leaves  are interchanged, and the local
\(L/R\)-labels are transported along the reversed edges. Clearly,
\(\operatorname{Rev}(N)\in\mathcal N_{n,\ell,k}\), and
\(\operatorname{Rev}\) is an involution on $\mathcal N_{n}:=\bigcup\limits_{0\leq k,\ell\leq n}\mathcal N_{n,k,\ell}$.
Now set
\[
\Phi
=
\Psi^{-1}
\circ
\operatorname{Rev}
\circ
\Psi,
\]
which proves Theorem~\ref{thm:bijective-symmetry}. 
\end{proof}

It remains to construct the bijection $\Psi$. To this end, we need the following canonical edge postorder for 
decorated binary networks. 
\begin{definition}
\label{def:postorder}
Let $N$ be a decorated binary network with $M+1$ edges. Suppose that $N$ has  $q$ roots $r_1,\ldots,r_q$, listed  in increasing order of their labels. 

We construct an ordered sequence $E$ of edges recursively. Initially, all vertices are unmarked and the sequence is empty. To process an edge $e=(v,u)$, if $u$ is unmarked, mark $u$ and recursively process its outgoing edges in the prescribed order ($L$ before $R$ at a tree vertex). 
Afterwards,  append $e$ to the sequence $E$.

Begin by processing the unique outgoing edge of $r_1$. Once this recursive process is completed, proceed to the unique outgoing edge of $r_2$, and continue in this way through $r_q$, retaining all marks throughout.
The resulting sequence
\[
E=(e_0,e_1,\ldots,e_M)
\]
is called the \blue{\it canonical edge postorder} of $N$.
\end{definition}

\begin{example}\label{ex:inverse-523}
Fig.~\ref{fig:inverse-523-postorder} shows a decorated binary network
$N\in\mathcal N_{5,2,3}$ with its canonical edge postorder
$e_0,e_1,\ldots,e_{10}$. In particular, recursively processing the outgoing edges of $r_1,r_2,r_3,r_4$ yields the segments $e_0$ to $e_2$, $e_3$ to $e_8$, $e_9$ and $e_{10}$, respectively. 
\end{example}

\begin{figure}
\noindent\makebox[\linewidth][c]{%
\begin{minipage}[c]{5.85cm}
\centering
\begin{tikzpicture}[
    x=1.10cm, y=0.79cm,
    edge/.style={draw=black,line width=0.49pt,line cap=round},
    port/.style={font=\fontsize{8}{9}\selectfont,inner sep=0pt,outer sep=0pt},
    boundary/.style={font=\fontsize{8.8}{9.8}\selectfont,inner sep=0pt,outer sep=0pt},
    edge label/.style={font=\fontsize{8.5}{9.5}\selectfont,inner sep=0.5pt,outer sep=0pt}
]

  % Vertices
  \coordinate (z1) at (0,0.50);
  \coordinate (z2) at (0,1.55);
  \coordinate (z3) at (1.10,3.50);
  \coordinate (t1) at (1.10,2.50);
  \coordinate (t2) at (1.10,4.90);
  \coordinate (r1) at (-1.08,2.35);
  \coordinate (r2) at (1.10,5.90);
  \coordinate (r3) at (-1.30,1.20);
  \coordinate (r4) at (3.10,1.55);
  \coordinate (l1) at (3.10,0.55);
  \coordinate (l2) at (0,-0.50);
  \coordinate (l3) at (2.12,1.55);

  % Edges
  \draw[edge] (z1)--(l2);
  \draw[edge] (z2)--(z1);
  \draw[edge] (r1)--(z2);
  \draw[edge] (t1)--(z2);
  \draw[edge] (t1)--(l3);
  \draw[edge] (z3)--(t1);
  \draw[edge] (t2)
    .. controls (0.45,4.58) and (0.45,3.82) .. (z3);
  \draw[edge] (t2)
    .. controls (1.75,4.58) and (1.75,3.82) .. (z3);
  \draw[edge] (r2)--(t2);
  \draw[edge] (r3)--(z1);
  \draw[edge] (r4)--(l1);

  % Vertices and boundary labels
  \foreach \v in {z1,z2,z3,t1,t2,r1,r2,r3,r4,l1,l2,l3}{
    \fill[black] (\v) circle[radius=2.30pt];
  }

  \foreach \v/\i in {r1/1,r2/2,r3/3,r4/4}{
    \node[boundary,above=3.2pt] at (\v) {$r_{\i}$};
  }

  \foreach \v/\i in {l1/1,l2/2,l3/3}{
    \node[boundary,text=wallred,below=3.2pt] at (\v) {$\i$};
  }

  % Local L/R labels
  \node[port] at (0.13,0.85) {$R$};
  \node[port] at (-0.32,0.5) {$L$};
  \node[port] at (-0.35,1.62) {$R$};
  \node[port] at (0.34,1.67) {$L$};
  
  \node[port] at (0.81,2.44) {$L$};
  \node[port] at (1.41,2.44) {$R$};
  
  \node[port] at (0.69,4.78) {$L$};
  \node[port] at (1.53,4.78) {$R$};
  \node[port] at (0.7,3.65) {$L$};
  \node[port] at (1.5,3.65) {$R$};

  % Canonical edge postorder
  \magenta{\node[edge label,anchor=east] at (0,0) {$e_0$};
  \node[edge label,anchor=east] at (0,1.06) {$e_1$};
  \node[edge label] at (-0.69,2.24) {$e_2$};
  \node[edge label] at (0.38,2.12) {$e_3$};
  \node[edge label] at (1.76,2.12) {$e_4$};
  \node[edge label,anchor=east] at (1.1,3.00) {$e_5$};
  \node[edge label,anchor=east] at (0.6,4.20) {$e_6$};
  \node[edge label,anchor=west] at (1.6,4.20) {$e_7$};
  \node[edge label,anchor=east] at (1.1,5.40) {$e_8$};
  \node[edge label] at (-0.85,1.12) {$e_9$};
  \node[edge label,anchor=west] at (3.12,1.05) {$e_{10}$};
}
\end{tikzpicture}
\end{minipage}%
\hspace{0.65cm}%
\begin{minipage}[c]{4.25cm}
\centering

{\renewcommand{\arraystretch}{1.60}
\begin{tabular*}{4.25cm}{@{\extracolsep{\fill}}ll@{}}
  \toprule
  \textbf{Type} & \textbf{Edges}\\
  \midrule
  $\mathsf A$ & $e_4,\ e_{10}$\\
  $\mathsf B$ & $e_1,\ e_2,\ e_6$\\
  $\mathsf C$ & $e_3,\ e_7,\ e_9$\\
  $\mathsf D$ & $e_5,\ e_8$\\
  \bottomrule
\end{tabular*}\par}
\end{minipage}
}
\captionsetup{hypcap=false}
\caption{The canonical edge postorder and edge types of  a decorated binary network. \label{fig:inverse-523-postorder}}
\end{figure}

We shall construct $\Psi$  by scanning the tableau
entries in increasing order and recover its inverse $\Psi^{-1}$
using the canonical edge postorder of decorated binary networks.

%%%%%%%%%%%%%%%%%%%%%%%
\subsection{The construction of $\Psi$}%%%
\label{subsec:increasing-scan}
%%%%%%%%%%%%%%%%%%%%

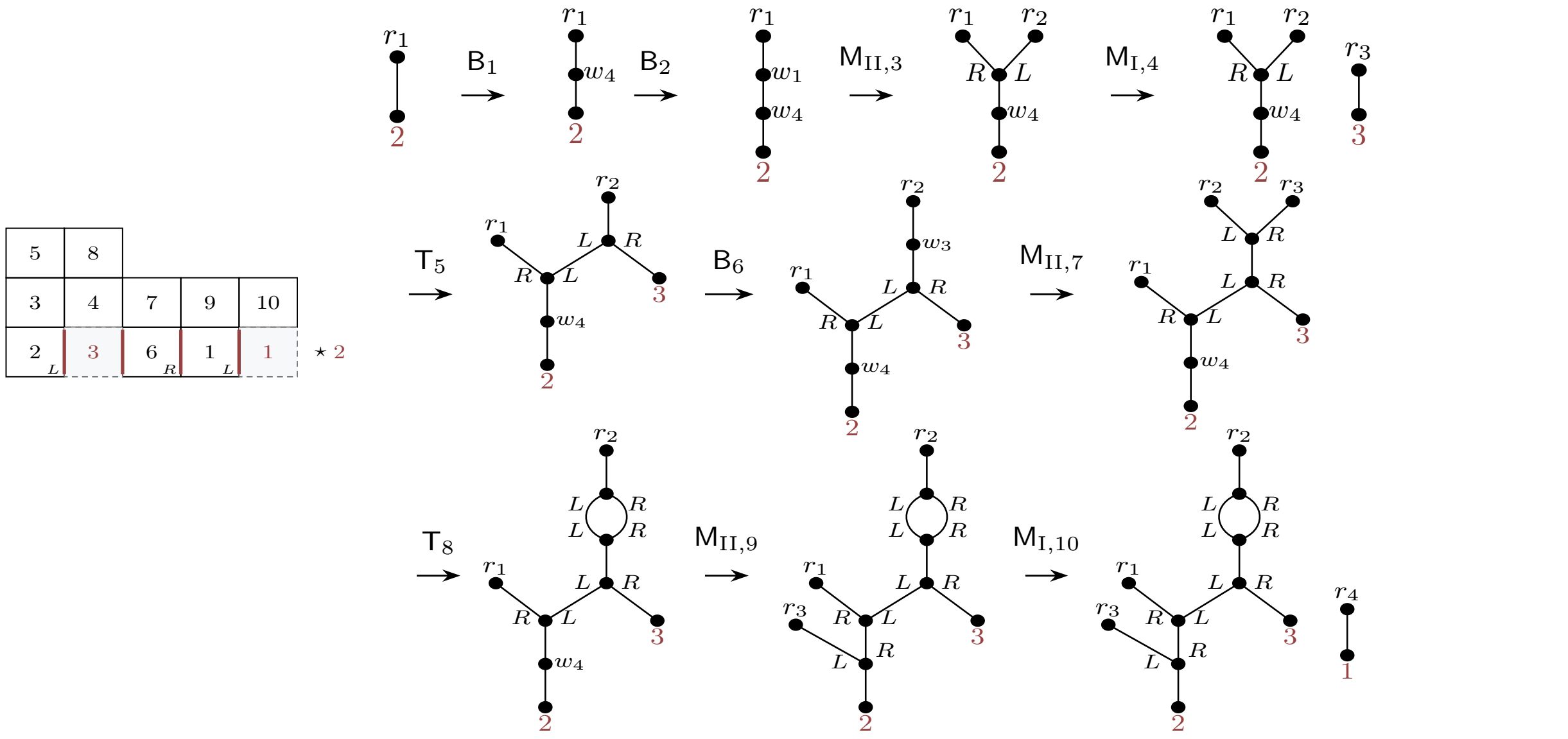
\begin{figure}%[h!]
\begin{tikzpicture}[yscale=0.85,
    x=1cm,y=1cm,
    dbn edge/.style={draw=black,line width=0.49pt,line cap=round},
    dbn port/.style={font=\fontsize{7.5}{8.1}\selectfont,
      inner sep=0pt,outer sep=0pt},
    dbn boundary/.style={font=\fontsize{8.8}{9.8}\selectfont,
      inner sep=0pt,outer sep=0pt},
    dbn temp/.style={font=\fontsize{7.7}{8.5}\selectfont,
      inner sep=0pt,outer sep=0pt},
    fwd arrow/.style={-{Stealth[length=1.85mm,width=1.25mm]},line width=0.54pt}
]
\path[use as bounding box] (0.10,1.58) rectangle (14.44,10.52);

\def\ForwardTableauScale{0.75}
\begin{scope}[shift={(0.5cm,6cm)}]
  \begin{scope}[
    % Physical ratio W:H ~ 1.16:1 after the outer yscale=0.85.
    % Preserve one identical geometry for every cell, including bottom L/R.
    x=0.8cm,y=0.8cm,
    scale=\ForwardTableauScale,transform shape,
    every node/.style={font=\fontsize{9}{10}\selectfont,
                       inner sep=0.5pt,outer sep=0pt},
    tableau cell/.style={draw=black,line width=0.37pt},
    hole/.style={dash pattern=on 1.74pt off 1.22pt,
                 draw=holegray,fill=holefill,line width=0.33pt},
    wall/.style={draw=wallred,line width=1pt}
  ]
    % All rows use equal cell dimensions; aspect ratio is close to square.
    \foreach \c/\val in {1/5,2/8}{
      \draw[tableau cell] (\c-1,2) rectangle (\c,3);
      \node at (\c-.5,2.5) {$\val$};
    }
    \foreach \c/\val in {1/3,2/4,3/7,4/9,5/10}{
      \draw[tableau cell] (\c-1,1) rectangle (\c,2);
      \node at (\c-.5,1.5) {$\val$};
    }

    % Bottom row: exact same cell size as above.  L/R stays at bottom right.
    \foreach \c/\val/\eps in {1/2/L,3/6/R,4/1/L}{
      \draw[tableau cell] (\c-1,0) rectangle (\c,1);
      \node at (\c-.5,.5) {$\val$};
      \node[anchor=south east,font=\fontsize{6}{6.5}\selectfont]
        at (\c-.055,.035) {$\eps$};
    }
    \foreach \c/\val in {2/3,5/1}{
      \draw[hole] (\c-1,0) rectangle (\c,1);
      \node[text=wallred] at (\c-.5,.5) {$\val$};
    }
    \foreach \xx in {1,2,3,4}{
      \draw[wall] (\xx,.05)--(\xx,.95);
    }
\node at (5.4,.5) {$\star$};
\node[text=wallred, anchor=west] at (5.6,.5) {$2$};
 
  \end{scope}
\end{scope}

\begin{scope}[shift={(3.9cm,0cm)}]

% Row 1: stages 0-4 (on the right)
% N_0
\begin{scope}[shift={(0.630cm,9.300cm)},x=0.72cm,y=0.72cm]
  \coordinate (r1) at (0,0.8);
  \coordinate (q2) at (0,-0.2);

  \draw[dbn edge] (r1)--(q2);

  \foreach \v in {r1,q2}{
    \fill[black] (\v) circle[radius=2.30pt];
  }
  \node[dbn boundary,above=3pt] at (r1) {$r_1$};
  \node[dbn boundary,text=wallred,below=3.2pt]
    at (q2) {$2$};

  % Decorations: each label stays adjacent to its incident arc.

\end{scope}
\draw[fwd arrow] (1.290,9.410) -- (1.740,9.410)
  node[midway,above=2.1pt,font=\scriptsize] {$\mathsf B_1$};
% N_1
\begin{scope}[shift={(2.470cm,8.980cm)},x=0.72cm,y=0.72cm]
  \coordinate (r1) at (0,1.60);
  \coordinate (w4) at (0,0.95);
  \coordinate (q2) at (0,0.30);

  \draw[dbn edge] (r1)--(w4)--(q2);

  \foreach \v in {r1,w4,q2}{
    \fill[black] (\v) circle[radius=2.30pt];
  }
  \node[dbn boundary,above=3.2pt] at (r1) {$r_1$};
  \node[dbn boundary,text=wallred,below=3.2pt]
    at (q2) {$2$};
  \node[dbn temp,anchor=west] at (0.11,0.95) {$w_4$};

  % Decorations: each label stays adjacent to its incident arc.

\end{scope}
\draw[fwd arrow] (3.070,9.410) -- (3.520,9.410)
  node[midway,above=2.1pt,font=\scriptsize] {$\mathsf B_2$};
% N_2
\begin{scope}[shift={(4.400cm,8.580cm)},x=0.72cm,y=0.72cm]
  \coordinate (r1) at (0,2.15);
  \coordinate (w1) at (0,1.50);
  \coordinate (w4) at (0,0.85);
  \coordinate (q2) at (0,0.20);

  \draw[dbn edge] (r1)--(w1)--(w4)--(q2);

  \foreach \v in {r1,w1,w4,q2}{
    \fill[black] (\v) circle[radius=2.30pt];
  }
  \node[dbn boundary,above=3.2pt] at (r1) {$r_1$};
  \node[dbn boundary,text=wallred,below=3.2pt]
    at (q2) {$2$};

  \node[dbn temp,anchor=west] at (0.11,1.50) {$w_1$};
  \node[dbn temp,anchor=west] at (0.11,0.85) {$w_4$};

  % Decorations: each label stays adjacent to its incident arc.

\end{scope}
\draw[fwd arrow] (5.290,9.410) -- (5.740,9.410)
  node[midway,above=2.1pt,font=\scriptsize] {$\mathsf M_{\mathrm{II},3}$};
% N_3
\begin{scope}[shift={(6.830cm,8.580cm)},x=0.72cm,y=0.72cm]
  \coordinate (r1) at (-0.52,2.15);
  \coordinate (r2) at ( 0.52,2.15);
  \coordinate (w1) at (0,1.50);
  \coordinate (w4) at (0,0.85);
  \coordinate (q2) at (0,0.20);

  \draw[dbn edge] (r1)--(w1);
  \draw[dbn edge] (r2)--(w1);
  \draw[dbn edge] (w1)--(w4)--(q2);

  \foreach \v in {r1,r2,w1,w4,q2}{
    \fill[black] (\v) circle[radius=2.30pt];
  }
  \foreach \v/\i in {r1/1,r2/2}{
    \node[dbn boundary,above=3.2pt] at (\v) {$r_{\i}$};
  }
  \node[dbn boundary,text=wallred,below=3.2pt]
    at (q2) {$2$};

  \node[dbn temp,anchor=west] at (0.11,0.85) {$w_4$};

  % Decorations: each label stays adjacent to its incident arc.
  \node[dbn port] at ($(w1)+(-0.34,0.04)$) {$R$};
  \node[dbn port] at ($(w1)+(0.34,0.04)$) {$L$};
\end{scope}
\draw[fwd arrow] (7.980,9.410) -- (8.430,9.410)
  node[midway,above=2.1pt,font=\scriptsize] {$\mathsf M_{\mathrm I,4}$};
% N_4
\begin{scope}[shift={(9.600cm,8.580cm)},x=0.72cm,y=0.72cm]
  \coordinate (r1) at (-0.62,2.15);
  \coordinate (r2) at ( 0.42,2.15);
  \coordinate (w1) at (-0.10,1.50);
  \coordinate (w4) at (-0.10,0.85);
  \coordinate (q2) at (-0.10,0.20);
  \coordinate (r3) at (1.30,1.58);
  \coordinate (q3) at (1.30,0.83);

  \draw[dbn edge] (r1)--(w1);
  \draw[dbn edge] (r2)--(w1);
  \draw[dbn edge] (w1)--(w4)--(q2);
  \draw[dbn edge] (r3)--(q3);

  \foreach \v in {r1,r2,r3,w1,w4,q2,q3}{
    \fill[black] (\v) circle[radius=2.30pt];
  }
  \foreach \v/\i in {r1/1,r2/2,r3/3}{
    \node[dbn boundary,above=3.2pt] at (\v) {$r_{\i}$};
  }
  \foreach \v/\i in {q2/2,q3/3}{
    \node[dbn boundary,text=wallred,below=3.2pt]
      at (\v) {$\i$};
  }

  \node[dbn temp,anchor=west] at (0.01,0.85) {$w_4$};

  % Decorations: each label stays adjacent to its incident arc.
  \node[dbn port] at ($(w1)+(-0.34,0.04)$) {$R$};
  \node[dbn port] at ($(w1)+(0.34,0.04)$) {$L$};
\end{scope}

% The size of the second and third network rows is adjustable here.
% First network row stays at full size
\def\ForwardLowerRowsScale{0.92}
% Row 2: stages 5-7 (uniform network scaling; shifted up by 0.30 cm)
\draw[fwd arrow] (0.75,7) -- (1.2,7)
  node[midway,above=2.1pt,font=\scriptsize] {$\mathsf T_5$};
% N_5
\begin{scope}[shift={(2.420cm,6.250cm)},x=0.76cm,y=0.76cm,scale=\ForwardLowerRowsScale,transform shape]
  \coordinate (w1) at (-0.35,1.35);
  \coordinate (w4) at (-0.35,0.60);
  \coordinate (q2) at (-0.35,-0.15);
  \coordinate (v1) at ( 0.55,2.00);
  \coordinate (r2) at ( 0.55,2.75);
  \coordinate (r1) at (-1.08,2.00);
  \coordinate (q3) at ( 1.30,1.35);

  \draw[dbn edge] (r1)--(w1);
  \draw[dbn edge] (r2)--(v1);
  \draw[dbn edge] (v1)--(w1);
  \draw[dbn edge] (v1)--(q3);
  \draw[dbn edge] (w1)--(w4)--(q2);

  \foreach \v in {r1,r2,v1,w1,w4,q2,q3}{
    \fill[black] (\v) circle[radius=2.30pt];
  }
  \foreach \v/\i in {r1/1,r2/2}{
    \node[dbn boundary,above=3.2pt] at (\v) {$r_{\i}$};
  }
  \foreach \v/\i in {q2/2,q3/3}{
    \node[dbn boundary,text=wallred,below=3.2pt]
      at (\v) {$\i$};
  }

  \node[dbn temp,anchor=west] at (-0.22,0.60) {$w_4$};

  % Decorations: each label stays adjacent to its incident arc.
  \node[dbn port] at ($(w1)+(-0.36,0.06)$) {$R$};
  \node[dbn port] at ($(w1)+(0.34,0.06)$) {$L$};
  \node[dbn port] at ($(v1)+(-0.35,0.02)$) {$L$};
  \node[dbn port] at ($(v1)+(0.35,0.02)$) {$R$};
\end{scope}
\draw[fwd arrow] (3.8,7) -- (4.3,7)
  node[midway,above=2.1pt,font=\scriptsize] {$\mathsf B_6$};
% N_6
\begin{scope}[shift={(5.560cm,5.680cm)},x=0.76cm,y=0.76cm,scale=\ForwardLowerRowsScale,transform shape]
  \coordinate (w1) at (-0.35,1.35);
  \coordinate (w4) at (-0.35,0.60);
  \coordinate (q2) at (-0.35,-0.15);
  \coordinate (v1) at ( 0.55,2.00);
  \coordinate (w3) at ( 0.55,2.75);
  \coordinate (r2) at ( 0.55,3.50);
  \coordinate (r1) at (-1.08,2.00);
  \coordinate (q3) at ( 1.30,1.35);

  \draw[dbn edge] (r1)--(w1);
  \draw[dbn edge] (r2)--(w3)--(v1);
  \draw[dbn edge] (v1)--(w1);
  \draw[dbn edge] (v1)--(q3);
  \draw[dbn edge] (w1)--(w4)--(q2);

  \foreach \v in {r1,r2,w3,v1,w1,w4,q2,q3}{
    \fill[black] (\v) circle[radius=2.30pt];
  }
  \foreach \v/\i in {r1/1,r2/2}{
    \node[dbn boundary,above=3.2pt] at (\v) {$r_{\i}$};
  }
  \foreach \v/\i in {q2/2,q3/3}{
    \node[dbn boundary,text=wallred,below=3.2pt]
      at (\v) {$\i$};
  }

  \node[dbn temp,anchor=west] at (0.68,2.75) {$w_3$};
  \node[dbn temp,anchor=west] at (-0.22,0.60) {$w_4$};

  % Decorations: each label stays adjacent to its incident arc.
  \node[dbn port] at ($(w1)+(-0.36,0.06)$) {$R$};
  \node[dbn port] at ($(w1)+(0.34,0.06)$) {$L$};
  \node[dbn port] at ($(v1)+(-0.35,0.02)$) {$L$};
  \node[dbn port] at ($(v1)+(0.35,0.02)$) {$R$};
\end{scope}
\draw[fwd arrow] (7.15,7) -- (7.6,7)
  node[midway,above=2.1pt,font=\scriptsize] {$\mathsf M_{\mathrm{II},7}$};
% N_7
\begin{scope}[shift={(9.050cm,5.750cm)},x=0.76cm,y=0.76cm,scale=\ForwardLowerRowsScale,transform shape]
  \coordinate (w1) at (-0.35,1.35);
  \coordinate (w4) at (-0.35,0.60);
  \coordinate (q2) at (-0.35,-0.15);
  \coordinate (v1) at ( 0.55,2.00);
  \coordinate (w3) at ( 0.55,2.75);
  \coordinate (r2) at (-0.05,3.40);
  \coordinate (r3) at ( 1.15,3.40);
  \coordinate (r1) at (-1.08,2.00);
  \coordinate (q3) at ( 1.30,1.35);

  \draw[dbn edge] (r1)--(w1);
  \draw[dbn edge] (r2)--(w3);
  \draw[dbn edge] (r3)--(w3);
  \draw[dbn edge] (w3)--(v1);
  \draw[dbn edge] (v1)--(w1);
  \draw[dbn edge] (v1)--(q3);
  \draw[dbn edge] (w1)--(w4)--(q2);

  \foreach \v in {r1,r2,r3,w3,v1,w1,w4,q2,q3}{
    \fill[black] (\v) circle[radius=2.30pt];
  }
  \foreach \v/\i in {r1/1,r2/2,r3/3}{
    \node[dbn boundary,above=3.2pt] at (\v) {$r_{\i}$};
  }
  \foreach \v/\i in {q2/2,q3/3}{
    \node[dbn boundary,text=wallred,below=3.2pt]
      at (\v) {$\i$};
  }

  \node[dbn temp,anchor=west] at (-0.22,0.60) {$w_4$};

  % Decorations: each label stays adjacent to its incident arc.
  \node[dbn port] at ($(w1)+(-0.36,0.06)$) {$R$};
  \node[dbn port] at ($(w1)+(0.34,0.06)$) {$L$};
  \node[dbn port] at ($(v1)+(-0.35,0.02)$) {$L$};
  \node[dbn port] at ($(v1)+(0.35,0.02)$) {$R$};
  \node[dbn port] at ($(w3)+(-0.34,0.08)$) {$L$};
  \node[dbn port] at ($(w3)+(0.34,0.08)$) {$R$};
\end{scope}

% Row 3: stages 8-10 (scale 0.92; shifted up by 0.60 cm)
\draw[fwd arrow] (0.830,3.590) -- (1.280,3.590)
  node[midway,above=2.1pt,font=\scriptsize] {$\mathsf T_8$};
% N_8
\begin{scope}[shift={(2.400cm,2.100cm)},x=0.76cm,y=0.76cm,scale=\ForwardLowerRowsScale,transform shape]
  \coordinate (w1) at (-0.35,1.35);
  \coordinate (w4) at (-0.35,0.60);
  \coordinate (q2) at (-0.35,-0.15);
  \coordinate (v1) at ( 0.55,2.00);
  \coordinate (w3) at ( 0.55,2.75);
  \coordinate (v2) at ( 0.55,3.55);
  \coordinate (r2) at ( 0.55,4.30);
  \coordinate (r1) at (-1.08,2.00);
  \coordinate (q3) at ( 1.30,1.35);

  \draw[dbn edge] (r1)--(w1);
  \draw[dbn edge] (r2)--(v2);

  % Parallel edges
  \draw[dbn edge] (v2)
    .. controls (0.15,3.40) and (0.15,2.90) .. (w3);
  \draw[dbn edge] (v2)
    .. controls (0.95,3.40) and (0.95,2.90) .. (w3);

  \draw[dbn edge] (w3)--(v1);
  \draw[dbn edge] (v1)--(w1);
  \draw[dbn edge] (v1)--(q3);
  \draw[dbn edge] (w1)--(w4)--(q2);

  \foreach \v in {r1,r2,v2,w3,v1,w1,w4,q2,q3}{
    \fill[black] (\v) circle[radius=2.30pt];
  }
  \foreach \v/\i in {r1/1,r2/2}{
    \node[dbn boundary,above=3.2pt] at (\v) {$r_{\i}$};
  }
  \foreach \v/\i in {q2/2,q3/3}{
    \node[dbn boundary,text=wallred,below=3.2pt]
      at (\v) {$\i$};
  }

  \node[dbn temp,anchor=west] at (-0.22,0.60) {$w_4$};

  % Decorations: each label stays adjacent to its incident arc.
  \node[dbn port] at ($(w1)+(-0.36,0.06)$) {$R$};
  \node[dbn port] at ($(w1)+(0.34,0.06)$) {$L$};
  \node[dbn port] at ($(v1)+(-0.35,0.02)$) {$L$};
  \node[dbn port] at ($(v1)+(0.35,0.02)$) {$R$};
  \node[dbn port] at ($(v2)+(-0.45,-0.17)$) {$L$};
  \node[dbn port] at ($(v2)+(0.45,-0.17)$) {$R$};
  \node[dbn port] at ($(w3)+(-0.45,0.17)$) {$L$};
  \node[dbn port] at ($(w3)+(0.45,0.17)$) {$R$};
\end{scope}
\draw[fwd arrow] (3.8,3.590) -- (4.25,3.590)
  node[midway,above=2.1pt,font=\scriptsize] {$\mathsf M_{\mathrm{II},9}$};
% N_9
\begin{scope}[shift={(5.700cm,2.100cm)},x=0.76cm,y=0.76cm,scale=\ForwardLowerRowsScale,transform shape]
  \coordinate (w1) at (-0.35,1.35);
  \coordinate (w4) at (-0.35,0.60);
  \coordinate (q2) at (-0.35,-0.15);
  \coordinate (v1) at ( 0.55,2.00);
  \coordinate (w3) at ( 0.55,2.75);
  \coordinate (v2) at ( 0.55,3.55);
  \coordinate (r2) at ( 0.55,4.30);
  \coordinate (r1) at (-1.08,2.00);
  \coordinate (r3) at (-1.38,1.28);
  \coordinate (q3) at ( 1.30,1.35);

  \draw[dbn edge] (r1)--(w1);
  \draw[dbn edge] (r2)--(v2);

  \draw[dbn edge] (v2)
    .. controls (0.15,3.40) and (0.15,2.90) .. (w3);
  \draw[dbn edge] (v2)
    .. controls (0.95,3.40) and (0.95,2.90) .. (w3);

  \draw[dbn edge] (w3)--(v1);
  \draw[dbn edge] (v1)--(w1);
  \draw[dbn edge] (v1)--(q3);

  \draw[dbn edge] (w1)--(w4);
  \draw[dbn edge] (r3)--(w4);
  \draw[dbn edge] (w4)--(q2);

  \foreach \v in {r1,r2,r3,v2,w3,v1,w1,w4,q2,q3}{
    \fill[black] (\v) circle[radius=2.30pt];
  }
  \foreach \v/\i in {r1/1,r2/2,r3/3}{
    \node[dbn boundary,above=3.2pt] at (\v) {$r_{\i}$};
  }
  \foreach \v/\i in {q2/2,q3/3}{
    \node[dbn boundary,text=wallred,below=3.2pt]
      at (\v) {$\i$};
  }

  % Decorations: each label stays adjacent to its incident arc.
  \node[dbn port] at ($(w1)+(-0.36,0.06)$) {$R$};
  \node[dbn port] at ($(w1)+(0.34,0.06)$) {$L$};
  \node[dbn port] at ($(v1)+(-0.35,0.02)$) {$L$};
  \node[dbn port] at ($(v1)+(0.35,0.02)$) {$R$};
  \node[dbn port] at ($(v2)+(-0.45,-0.17)$) {$L$};
  \node[dbn port] at ($(v2)+(0.45,-0.17)$) {$R$};
  \node[dbn port] at ($(w3)+(-0.45,0.17)$) {$L$};
  \node[dbn port] at ($(w3)+(0.45,0.17)$) {$R$};
  \node[dbn port] at ($(w4)+(-0.39,0.04)$) {$L$};
  \node[dbn port] at ($(w4)+(0.28,0.24)$) {$R$};
\end{scope}
\draw[fwd arrow] (7.1,3.590) -- (7.54,3.590)
  node[midway,above=2.1pt,font=\scriptsize] {$\mathsf M_{\mathrm I,10}$};
% N_10
\begin{scope}[shift={(8.920cm,2.100cm)},x=0.76cm,y=0.76cm,scale=\ForwardLowerRowsScale,transform shape]
  \coordinate (w1) at (-0.35,1.35);
  \coordinate (w4) at (-0.35,0.60);
  \coordinate (q2) at (-0.35,-0.15);
  \coordinate (v1) at ( 0.55,2.00);
  \coordinate (w3) at ( 0.55,2.75);
  \coordinate (v2) at ( 0.55,3.55);
  \coordinate (r2) at ( 0.55,4.30);
  \coordinate (r1) at (-1.08,2.00);
  \coordinate (r3) at (-1.38,1.28);
  \coordinate (q3) at ( 1.30,1.35);
  \coordinate (r4) at ( 2.14,1.55);
  \coordinate (q1) at ( 2.14,0.75);

  \draw[dbn edge] (r1)--(w1);
  \draw[dbn edge] (r2)--(v2);

  \draw[dbn edge] (v2)
    .. controls (0.15,3.40) and (0.15,2.90) .. (w3);
  \draw[dbn edge] (v2)
    .. controls (0.95,3.40) and (0.95,2.90) .. (w3);

  \draw[dbn edge] (w3)--(v1);
  \draw[dbn edge] (v1)--(w1);
  \draw[dbn edge] (v1)--(q3);

  \draw[dbn edge] (w1)--(w4);
  \draw[dbn edge] (r3)--(w4);
  \draw[dbn edge] (w4)--(q2);
  \draw[dbn edge] (r4)--(q1);

  \foreach \v in {r1,r2,r3,r4,v2,w3,v1,w1,w4,q1,q2,q3}{
    \fill[black] (\v) circle[radius=2.30pt];
  }
  \foreach \v/\i in {r1/1,r2/2,r3/3,r4/4}{
    \node[dbn boundary,above=3.2pt] at (\v) {$r_{\i}$};
  }
  \foreach \v/\i in {q1/1,q2/2,q3/3}{
    \node[dbn boundary,text=wallred,below=3.2pt]
      at (\v) {$\i$};
  }

  % Decorations: each label stays adjacent to its incident arc.
  \node[dbn port] at ($(w1)+(-0.36,0.06)$) {$R$};
  \node[dbn port] at ($(w1)+(0.34,0.06)$) {$L$};
  \node[dbn port] at ($(v1)+(-0.35,0.02)$) {$L$};
  \node[dbn port] at ($(v1)+(0.35,0.02)$) {$R$};
  \node[dbn port] at ($(v2)+(-0.45,-0.17)$) {$L$};
  \node[dbn port] at ($(v2)+(0.45,-0.17)$) {$R$};
  \node[dbn port] at ($(w3)+(-0.45,0.17)$) {$L$};
  \node[dbn port] at ($(w3)+(0.45,0.17)$) {$R$};
  \node[dbn port] at ($(w4)+(-0.39,0.04)$) {$L$};
  \node[dbn port] at ($(w4)+(0.28,0.24)$) {$R$};
\end{scope}
\end{scope}
\end{tikzpicture}
\caption{An example of $\Psi$, where arrows labeled $\mathsf B_a$,
$\mathsf M_{\mathrm I,a}$, $\mathsf M_{\mathrm{II},a}$, and
$\mathsf T_a$ mark the bottom, Type~I middle, Type~II middle,
and top operations at entry $a$, respectively. \label{fig:Psi}}
\end{figure}

Fix \(X=(T;\varepsilon,\rho)\in\widehat{\mathcal D}_{n,k,\ell}\) and
set \(M=n+k+\ell\). We construct the intermediate digraphs 
\(N_0,N_1,\ldots,N_M\) as follows and finally set $\Psi(X)=N_M$. 
\begin{enumerate}[
    label=\textup{(\arabic*)},
    leftmargin=2.5em,
    itemsep=4pt
]
\item \textbf{Initialization.}
Set \(N_0\) to consist of the root \(r_1\), the leaf
$\rho(\star)$, and the edge
\(r_1\to\rho(\star)\). 
\item For \(1\le a\le M\), suppose that \(N_{a-1}\) has been constructed
and has roots \(r_1,\ldots,r_q\). We construct
\(N_a\) from \(N_{a-1}\) according to one of the following four cases.
\begin{enumerate}[
    leftmargin=3em,
    itemsep=3pt,
    topsep=3pt
]
\item[\textup{(\(\mathsf M_{\mathrm I}\))}]
\uline{\textbf{Type I middle entry}: \(a\) is the middle entry in column \(j\) with the hole
\(h_j\) below it.} Add the root \(r_{q+1}\), the leaf
\(\rho(h_j)\), and the edge
\(r_{q+1}\to \rho(h_j)\). 
See Fig.~\ref{fig:forward-cases}\subref{fig:forward-case-a} for the illustration. 

\item[\textup{(\(\mathsf B\))}]
\uline{\textbf{Bottom entry}:
 \(a\) lies in the occupied bottom cell of column \(j\).}
Let \(r_q\to u\) be the outgoing edge of \(r_q\). Replace this edge by
\(r_q\to w_j\to u\), where \(w_j\) is a new vertex. If \(u\) is a reticulation vertex, transfer the
local label at \(u\) by setting
\(\tau_u(w_j\to u)=\tau_u(r_q\to u)\).
Remember \(w_j\) until the middle entry in column \(j\) is processed.
See Fig.~\ref{fig:forward-cases}\subref{fig:forward-case-b} for the illustration.

\item[\textup{(\(\mathsf M_{\mathrm{II}}\))}]
\uline{\textbf{Type II middle entry}: \(a\) is the middle entry in column \(j\), whose bottom
cell is occupied.} Let \(e\) be the unique incoming edge of \(w_j\).
Add the root \(r_{q+1}\) and the edge \(r_{q+1}\to w_j\). Define
\(\tau_{w_j}(r_{q+1}\to w_j)=\varepsilon(j)\), and give \(e\) the
other label in \(\{L,R\}\).
See Fig.~\ref{fig:forward-cases}\subref{fig:forward-case-c} for the illustration.

\item[\textup{(\(\mathsf T\))}]
\uline{\textbf{Top entry}: \(a\) lies in a top cell.} Let
\(r_{q-1}\to x\) and \(r_q\to y\) be the outgoing edges of the last
two roots (possibly $x=y$). Delete these two roots and their outgoing edges, and replace
them by a new root, again denoted \(r_{q-1}\), a tree vertex \(v\),
and the edges \(r_{q-1}\to v\), \(v\to x\), and \(v\to y\).
Set \(\sigma_v(v\to x)=L\) and
\(\sigma_v(v\to y)=R\).
If \(x\) is a reticulation vertex, set
\(\tau_x(v\to x)=\tau_x(r_{q-1}\to x)\). Likewise, if \(y\) is a
reticulation vertex, set
\(\tau_y(v\to y)=\tau_y(r_q\to y)\).
See Fig.~\ref{fig:forward-cases}\subref{fig:forward-case-d} for the illustration.
\end{enumerate}

\begin{figure}
\begin{subfigure}[t]{0.225\textwidth}
\centering
\begin{tikzpicture}[
    x=1cm, y=1cm,
    line width=0.44pt,
    line cap=round,
    every node/.style={font=\small,inner sep=0.6pt},
    vertex/.style={circle,draw=black,fill=black,
        inner sep=0pt,outer sep=0pt,minimum size=3.1pt},
    change/.style={font=\normalsize}
]
  \node[vertex] (a0r) at (0.15,1.50) {};
  \node[vertex] (a0u) at (0.15,0.45) {};
  \draw (a0r)--(a0u);
  \node[above=2pt] at (a0r) {$r_q$};
  \node[below=2pt] at (a0u) {$u$};

  % Short arrow for the first case.
  \node[change] at (0.76,0.98) {$\rightarrow$};

  \node[vertex] (a1r) at (1.55,1.50) {};
  \node[vertex] (a1u) at (1.55,0.45) {};
  \node[vertex] (a1nr) at (2.42,1.50) {};
  \node[vertex] (a1leaf) at (2.42,0.45) {};
  \draw (a1r)--(a1u);
  \draw (a1nr)--(a1leaf);
  \node[above=2pt] at (a1r) {$r_q$};
  \node[below=2pt] at (a1u) {$u$};
  \node[above=2pt] at (a1nr) {$r_{q+1}$};
  \node[below=2pt,text=wallred] at (a1leaf) {$\rho(h_j)$};
\end{tikzpicture}
\caption{\(\mathsf M_{\mathrm I}\) entry.}
\label{fig:forward-case-a}
\end{subfigure}\hfill
\begin{subfigure}[t]{0.175\textwidth}
\centering
\begin{tikzpicture}[
    x=1cm, y=1cm,
    line width=0.44pt,
    line cap=round,
    every node/.style={font=\small,inner sep=0.6pt},
    vertex/.style={circle,draw=black,fill=black,
        inner sep=0pt,outer sep=0pt,minimum size=3.1pt},
    change/.style={font=\normalsize}
]
  \node[vertex] (b0r) at (0.15,1.50) {};
  \node[vertex] (b0u) at (0.15,0.45) {};
  \draw (b0r)--node[pos=0.80,right=2pt] {$\alpha$} (b0u);
  \node[above=2pt] at (b0r) {$r_q$};
  \node[below=2pt] at (b0u) {$u$};

  % Short arrow for the second case.
  \node[change] at (0.72,0.98) {$\rightarrow$};

  \node[vertex] (b1r) at (1.38,1.50) {};
  \node[vertex] (b1w) at (1.38,0.98) {};
  \node[vertex] (b1u) at (1.38,0.45) {};
  \draw (b1r)--(b1w);
  \draw (b1w)--node[pos=0.72,right=2pt] {$\alpha$} (b1u);
  \node[above=2pt] at (b1r) {$r_q$};
  \node[right=2pt] at (b1w) {$w_j$};
  \node[below=2pt] at (b1u) {$u$};
\end{tikzpicture}
\caption{\(\mathsf B\) entry.}
\label{fig:forward-case-b}
\end{subfigure}\hfill
\begin{subfigure}[t]{0.255\textwidth}
\centering
\begin{tikzpicture}[
    x=1cm, y=1cm,
    line width=0.44pt,
    line cap=round,
    every node/.style={font=\small,inner sep=0.6pt},
    vertex/.style={circle,draw=black,fill=black,
        inner sep=0pt,outer sep=0pt,minimum size=3.1pt},
    change/.style={font=\normalsize}
]
  \node[vertex] (c0z) at (0.00,1.50) {};
  \node[vertex] (c0w) at (0.00,0.98) {};
  \node[vertex] (c0u) at (0.00,0.45) {};
  \draw (c0z)--(c0w)--(c0u);
  \node[above=2pt] at (c0z) {$z$};
  \node[right=2pt] at (c0w) {$w_j$};
  \node[below=2pt] at (c0u) {$u$};

  \node[change] at (0.92,0.98) {$\rightarrow$};

  \node[vertex] (c1z) at (1.85,1.50) {};
  \node[vertex] (c1r) at (3.05,1.50) {};
  \node[vertex] (c1w) at (2.45,0.88) {};
  \node[vertex] (c1u) at (2.45,0.38) {};
  \draw (c1z)--(c1w);
  \draw (c1r)--(c1w);
  \draw (c1w)--(c1u);
  \node[anchor=east] at (2.1,1.15) {$\bar\varepsilon(j)$};
  \node[anchor=west] at (2.78,1.2) {$\varepsilon(j)$};
  \node[above=2pt] at (c1z) {$z$};
  \node[above=2pt] at (c1r) {$r_{q+1}$};
  \node[right=2pt] at (c1w) {$w_j$};
  \node[below=2pt] at (c1u) {$u$};
\end{tikzpicture}
\caption{\(\mathsf M_{\mathrm{II}}\) entry.}
\label{fig:forward-case-c}
\end{subfigure}\hfill
\begin{subfigure}[t]{0.3\textwidth}
\centering
\begin{tikzpicture}[
    x=0.9cm, y=0.9cm,
    line width=0.44pt,
    line cap=round,
    every node/.style={font=\small,inner sep=0.6pt},
    vertex/.style={circle,draw=black,fill=black,
        inner sep=0pt,outer sep=0pt,minimum size=3.1pt},
    change/.style={font=\normalsize}
]
  \node[vertex] (d0r1) at (0.00,1.58) {};
  \node[vertex] (d0r2) at (0.95,1.58) {};
  \node[vertex] (d0x)  at (0.00,0.28) {};
  \node[vertex] (d0y)  at (0.95,0.28) {};
  \draw (d0r1)--(d0x);
  \draw (d0r2)--(d0y);
  \node[above=2pt] at (d0r1) {$r_{q-1}$};
  \node[above=2pt] at (d0r2) {$r_q$};
  \node[below=2pt] at (d0x) {$x$};
  \node[below=2pt] at (d0y) {$y$};

   \node[change] at (1.88,0.94) {$\rightarrow$};

  \node[vertex] (d1r) at (3.45,1.66) {};
  \node[vertex] (d1v) at (3.45,1.08) {};
  \node[vertex] (d1x) at (2.50,0.10) {};
  \node[vertex] (d1y) at (4.40,0.10) {};
  \draw (d1r)--(d1v);
  \draw (d1v)
      --node[pos=0.24,left=4pt] {$L$}
       (d1x);
  \draw (d1v)
      --node[pos=0.24,right=4pt] {$R$}
      (d1y);
  \node[above=2pt] at (d1r) {$r_{q-1}$};
  \node[above right=1pt] at (d1v) {$v$};
  \node[below=2pt] at (d1x) {$x$};
  \node[below=2pt] at (d1y) {$y$};
\end{tikzpicture}
\caption{\(\mathsf T\) entry.}
\label{fig:forward-case-d}
\end{subfigure}

\caption{Four local modifications in the construction of
$\Psi$. Here $\bar{\varepsilon}(j)$ denotes the opposite of
$\varepsilon(j)$ in $\{L,R\}$. \label{fig:forward-cases}}
\end{figure}
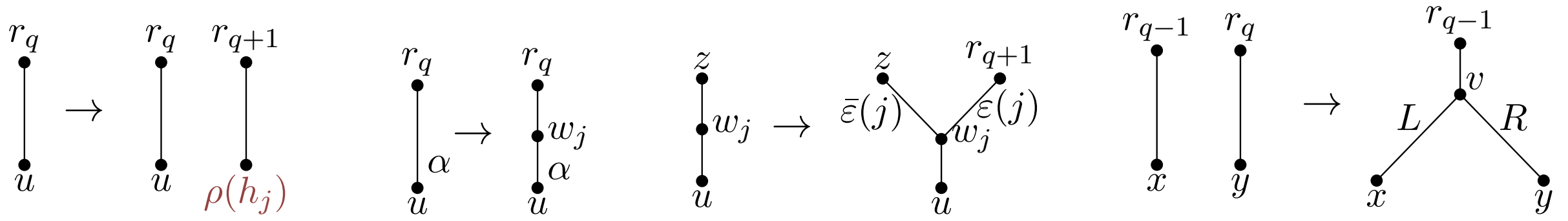
\item After \(N_M\) has been constructed, forget the temporary names of the
internal vertices.
The resulting decorated binary network is defined to be
\(\Psi(X)\).
\end{enumerate}
See Fig.~\ref{fig:Psi} for an example of $\Psi$. The interested reader is invited to consult  Appendices~\ref{sec:AppendixB} and~\ref{app:networks-211} for the correspondence  $\Psi:\widehat{\mathcal D}_{2,1,1}\longrightarrow\mathcal N_{2,1,1}$.

\begin{lemma}\label{prop:forward-well-defined}
The map
\(\Psi:\widehat{\mathcal D}_{n,k,\ell}\to\mathcal N_{n,k,\ell}\)
is well defined. Moreover, for any
\(X=(T;\varepsilon,\rho)\in\widehat{\mathcal D}_{n,k,\ell}\),
the network \(\Psi(X)\) has \(n-k+1\) roots and
\(n-\ell+1\) leaves.
\end{lemma}
\begin{proof}
To see that every step of $\Psi$ is well
defined, the only point requiring verification is Case~$\mathsf{T}$, where two roots
are needed. Suppose that the next entry \(a\) is a top entry, and that
\(m\) middle entries and \(t\) top entries have already been processed.
Since \(a\) itself is a top entry, the prefix \(1,\ldots,a\) contains
\(m\) middle entries and \(t+1\) top entries. By the tableau  conditions
on \(T\), we have \(m\ge t+1\).
It follows immediately from the construction that the current network
has \(m+1-t\) roots, since each middle step adds one root, each top
step removes one, and a bottom step preserves the number of roots. Since \(m+1-t\ge2\),  the two roots required in
Case~$\mathsf{T}$ are  available. Moreover, it is clear  from the construction that $\Psi(X)$ is acyclic,  and every vertex has an indegree–outdegree pair in
\(
\{(0,1), (1,0), (1,2), (2,1)\}.
\)
Thus, $\Psi(X)$ is a binary network and $\Psi$ is well defined.

At the end of the algorithm, since $T$ has \(n\) middle entries and \(k\) top entries,  $\Psi(X)$ has \(n-k+1\) roots.
Similarly, starting from one leaf and adding one precisely for each of
the \(n-\ell\) middle entries above a hole gives \(n-\ell+1\) leaves.
Thus the resulting network belongs to \(\mathcal N_{n,k,\ell}\).
\end{proof}

For an intermediate digraph \(N_a\) (\(0\le a\le M\)), we use the
same edge traversal as in Definition~\ref{def:postorder} and write
\[
E_a=(e_0^{(a)},e_1^{(a)},\ldots,e_a^{(a)})
\]
for the resulting canonical edge postorder.
For a vertex \(w\) of outdegree \(1\), order its incoming edges by
their occurrences in \(E_a\). If \(w\) has only one incoming edge,
we regard it as the earlier incoming edge; if \(w\) has two incoming edges, we call them the earlier and later incoming edges according to this order.

For \(1\le b\le a\), we assign to the edge  \(e_b^{(a)}=(v,u)\) of $E_a$ (excluding $e_0^{(a)}$) one of the following four types:
\[
\begin{array}{c|l}
\text{type} & \text{corresponding condition on } e_b^{(a)} \\ \hline
\mathsf A & u \text{ is a leaf}\\
\mathsf B & u \text{ has outdegree }1
              \text{ and }e_b^{(a)}\text{ is its earlier incoming edge}\\
\mathsf C & u \text{ has outdegree }1
              \text{ and }e_b^{(a)}\text{ is its later incoming edge}\\
\mathsf D & u \text{ is a tree vertex}
\end{array}
\]
For instance, the edge types of the network in Example~\ref{ex:inverse-523} are
shown in the right part of Fig.~\ref{fig:inverse-523-postorder}.

The following features of $\Psi$ are useful in constructing its inverse. 

\begin{lemma}\label{lem:edge-types}
With the notation above, the following assertions hold.
\begin{enumerate}
    \item[(i)]
    For \(1\le a\le M\), the last edge \(e_a^{(a)}\) of $E_a$ is of
    type \(\mathsf A\), \(\mathsf B\), \(\mathsf C\), or \(\mathsf D\) if and only if the entry \(a\) is a
    Type I middle entry, a bottom entry, a Type II middle entry,
    or a top entry in $T$, respectively.

    \item[(ii)]
    For \(1\le b\le a\le M\), the edges \(e_b^{(a)}\) and
    \(e_b^{(M)}\) have the same edge type.
\end{enumerate}
\end{lemma}

\begin{proof}
Part~(i) is straightforward from the construction.

 For~(ii), consider one step \(N_a\to N_{a+1}\). All edges not involved
  in the local modification keep their orders in the canonical
 edge postorder and their edge types.
It remains only to consider the edges involved
  in the local modification.  In Case~(\(\mathsf B\)),
 the edge \(r_q\to u\) is replaced by \(w_j\to u\); since the terminal
 vertex \(u\) is unchanged, the replacement edge remains the same type.
 In Case~($\mathsf{T}$), the edges \(r_{q-1}\to x\) and \(r_q\to y\) are replaced
 by \(v\to x\) and \(v\to y\), respectively; again, the terminal
 vertices are unchanged, so their edge types are preserved. Cases~(\(\mathsf M_{\mathrm{I}}\))
 and~(\(\mathsf M_{\mathrm{II}}\)) do not replace any previous edge.
Thus every previously occupied postorder position keeps the same edge
 type from one step to the next. Iterating this from \(N_a\) to \(N_M\)
 gives~(ii).
\end{proof}

 In the next section, we prove that $\Psi$ is a bijection by constructing its inverse explicitly. 
%%%%%%%%%%%%%%%%%%%%%%%%%%%%%%%%%%%%%%%%%%%%%%%%%%
\subsection{Constructing $\Psi^{-1}$ and finishing the proof of Theorem~\ref{thm:tableau-network-correspondence}}%%%%%%%%%%%%%%%%%%%%%%%%%%%%%%%%%%%%%%%
\label{subsec:inverse-map}
%%%%%%%%%%%%%%%%%%%%%%%%%%%%%%%%%%%%%%%%%%%%%%%%%%%

Fix \(N\in\mathcal N_{n,k,\ell}\), set \(M=n+k+\ell\), and let
\(E=(e_0,e_1,\ldots,e_M)\) be the canonical edge postorder of \(N\).
We construct \(\Psi^{-1}(N)\) by reading
\(e_1,e_2,\ldots,e_M\) in increasing order as follows.

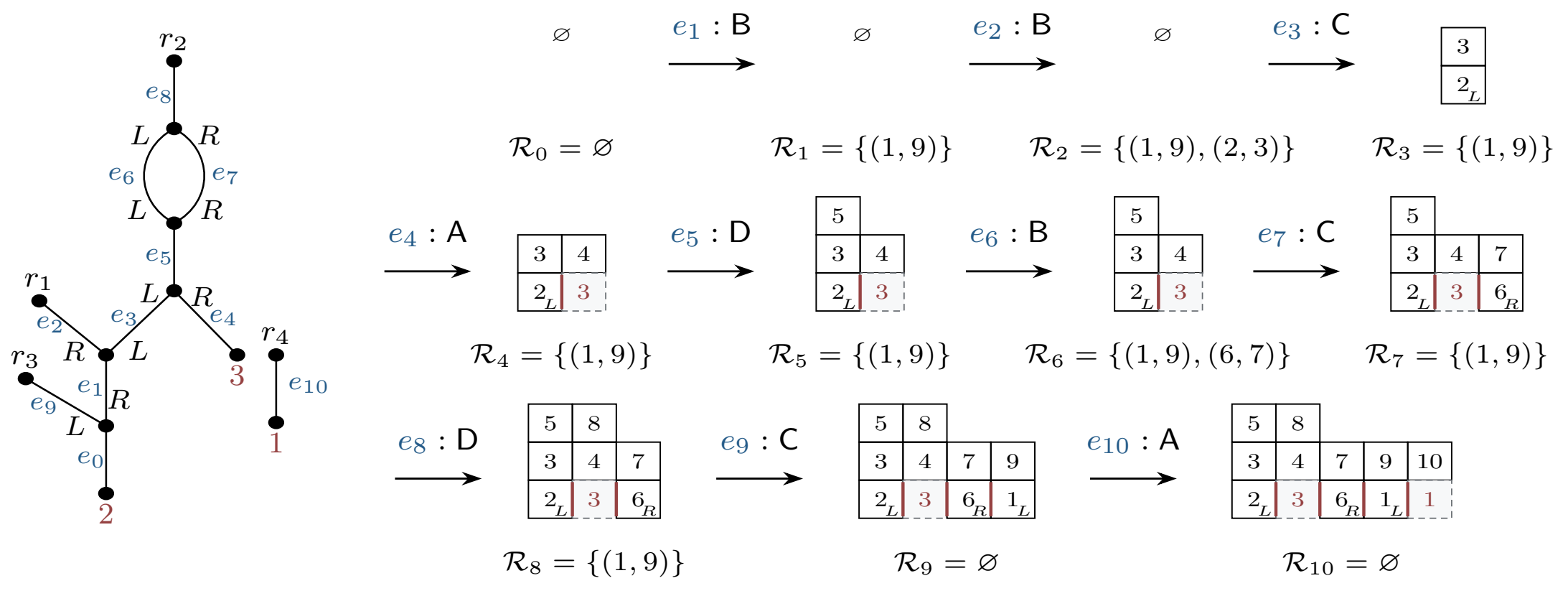
\begin{figure}
\centering
% MOD 1: Space between the drawing and its caption (local to this figure).
\setlength{\abovecaptionskip}{4pt}

\definecolor{edgeblue}{RGB}{36,91,137}
\def\InverseTableauScale{0.75}
% Center the whole drawing, including when it is wider than the text block.
\makebox[\linewidth][c]{%
\begin{tikzpicture}[yscale=0.83,
  x=1cm,y=1cm,
  cell/.style={draw=black,line width=0.37pt},
  hole/.style={dash pattern=on 1.74pt off 1.22pt,
               draw=holegray,fill=holefill,line width=0.33pt},
  wall/.style={draw=wallred,line width=0.83pt},
  stage/.style={font=\fontsize{8.5}{9.5}\selectfont,inner sep=0pt},
  recorded/.style={font=\fontsize{7.1}{8}\selectfont,inner sep=0pt},
  step arrow/.style={-{Stealth[length=1.85mm,width=1.25mm]},
                     line width=0.54pt}
]
% MOD 2: Use the natural TikZ bounding box; remove the forced bottom whitespace.
% Shift the input network further left; former x=0.80cm, new x=0.42cm.
\begin{scope}[shift={(1.7cm,4.55cm)},x=0.57cm,y=0.75cm,
  inv edge/.style={draw=black,line width=0.49pt,line cap=round},
  inv port/.style={font=\fontsize{7.2}{8}\selectfont,inner sep=0pt,outer sep=0pt},
  inv boundary/.style={font=\fontsize{7.8}{8.6}\selectfont,inner sep=0pt,outer sep=0pt},
  inv edge label/.style={font=\fontsize{7.0}{7.8}\selectfont,text=edgeblue,inner sep=0pt,outer sep=0pt}]
  \coordinate (z1) at (0,0.50);
  \coordinate (z2) at (0,1.55);
  \coordinate (z3) at (1.10,3.50);
  \coordinate (t1) at (1.10,2.50);
  \coordinate (t2) at (1.10,4.90);
  \coordinate (r1) at (-1.08,2.35);
  \coordinate (r2) at (1.10,5.90);
  \coordinate (r3) at (-1.30,1.20);
  \coordinate (r4) at (2.75,1.55);
  \coordinate (l1) at (2.75,0.55);
  \coordinate (l2) at (0,-0.50);
  \coordinate (l3) at (2.12,1.55);

  \draw[inv edge] (z1)--(l2);
  \draw[inv edge] (z2)--(z1);
  \draw[inv edge] (r1)--(z2);
  \draw[inv edge] (t1)--(z2);
  \draw[inv edge] (t1)--(l3);
  \draw[inv edge] (z3)--(t1);
  \draw[inv edge] (t2)
    .. controls (0.45,4.58) and (0.45,3.82) .. (z3);
  \draw[inv edge] (t2)
    .. controls (1.75,4.58) and (1.75,3.82) .. (z3);
  \draw[inv edge] (r2)--(t2);
  \draw[inv edge] (r3)--(z1);
  \draw[inv edge] (r4)--(l1);
  \foreach \v in {z1,z2,z3,t1,t2,r1,r2,r3,r4,l1,l2,l3}{
    \fill[black] (\v) circle[radius=2.10pt];
  }
  \foreach \v/\i in {r1/1,r2/2,r3/3,r4/4}{
    \node[inv boundary,above=2.8pt] at (\v) {$r_{\i}$};
  }
  \foreach \v/\i in {l1/1,l2/2,l3/3}{
    \node[inv boundary,text=wallred,below=2.8pt] at (\v) {$\i$};
  }
  % At binary vertices, port decorations retain the given L/R order.
  \node[inv port] at (0.2,0.90) {$R$};
  \node[inv port] at (-0.5,0.5) {$L$};
  \node[inv port] at (-0.53,1.61) {$R$};
  \node[inv port] at (0.52,1.62) {$L$};
  \node[inv port] at (0.7,2.47) {$L$};
  \node[inv port] at (1.54,2.41) {$R$};
  \node[inv port] at (0.55,4.8) {$L$};
  \node[inv port] at (1.65,4.8) {$R$};
  \node[inv port] at (0.5,3.7) {$L$};
  \node[inv port] at (1.7,3.7) {$R$};
  % Canonical edge postorder e_0,...,e_{10} (same as in the main text).
  \node[inv edge label,anchor=east] at (0,0.00) {$e_0$};
  \node[inv edge label,anchor=east] at (0,1.06) {$e_1$};
  \node[inv edge label] at (-0.9,2) {$e_2$};
  \node[inv edge label] at (0.3,2.1) {$e_3$};
  \node[inv edge label] at (1.9,2.1) {$e_4$};
  \node[inv edge label,anchor=east] at (1.1,3.00) {$e_5$};
  \node[inv edge label,anchor=east] at (0.5,4.20) {$e_6$};
  \node[inv edge label,anchor=west] at (1.7,4.20) {$e_7$};
  \node[inv edge label,anchor=east] at (1.1,5.40) {$e_8$};
  \node[inv edge label] at (-1,0.8) {$e_9$};
  \node[inv edge label,anchor=west] at (2.94,1.10) {$e_{10}$};
\end{scope}

% Shift the reconstructed tableaux and their arrows right as one unit.
% Changing this xshift does not affect the input network or Figure 5.
\begin{scope}[xshift=1.85cm]
% Three compact rows of tableaux, with the network at left:
\foreach \a/\panelx/\panely in {
   0/4.05/8.50, 1/6.82/8.50, 2/9.59/8.50, 3/12.36/8.50,
   4/4.05/6.20, 5/6.80/6.20, 6/9.55/6.20, 7/12.30/6.20,
   8/4.35/3.90, 9/7.60/3.90, 10/11.24/3.90
}{
  % Number of reconstructed middle cells at stage P_a.
  \def\ncols{0}
  \ifnum\a>2 \def\ncols{1}\fi
  \ifnum\a>3 \def\ncols{2}\fi
  \ifnum\a>6 \def\ncols{3}\fi
  \ifnum\a>8 \def\ncols{4}\fi
  \ifnum\a>9 \def\ncols{5}\fi

  % Each tableau is centered under its P_i and R_i labels.
  \pgfmathsetmacro{\stagehalf}{0.2025*\ncols}
  \pgfmathsetmacro{\stageleft}{\panelx-\stagehalf}
  \pgfmathsetmacro{\arrowhalf}{max(0.18,\stagehalf)}
  \coordinate (stageleft-\a) at ({\panelx-\arrowhalf},{\panely+0.44});
  \coordinate (stageright-\a) at ({\panelx+\arrowhalf},{\panely+0.44});

  \begin{scope}[shift={(\stageleft cm,\panely cm)}]
  \begin{scope}[
    x=0.54cm,y=0.563253cm,
    scale=\InverseTableauScale,transform shape,
    every node/.style={font=\fontsize{9}{10}\selectfont,
                       inner sep=0.5pt,outer sep=0pt}
  ]

    % Top row
    \ifnum\a>4
      \draw[cell] (0,2) rectangle (1,3);
      \node at (0.5,2.5) {$5$};
    \fi
    \ifnum\a>7
      \draw[cell] (1,2) rectangle (2,3);
      \node at (1.5,2.5) {$8$};
    \fi

    % Middle row
    \ifnum\a>2
      \draw[cell] (0,1) rectangle (1,2);
      \node at (0.5,1.5) {$3$};
    \fi
    \ifnum\a>3
      \draw[cell] (1,1) rectangle (2,2);
      \node at (1.5,1.5) {$4$};
    \fi
    \ifnum\a>6
      \draw[cell] (2,1) rectangle (3,2);
      \node at (2.5,1.5) {$7$};
    \fi
    \ifnum\a>8
      \draw[cell] (3,1) rectangle (4,2);
      \node at (3.5,1.5) {$9$};
    \fi
    \ifnum\a>9
      \draw[cell] (4,1) rectangle (5,2);
      \node at (4.5,1.5) {$10$};
    \fi

    % Bottom entries, decorations, and holes
    \ifnum\a>2
      \draw[cell] (0,0) rectangle (1,1);
      \node at (0.5,0.5) {$2$};
      \node[anchor=south east,
            font=\fontsize{6}{6.5}\selectfont]
        at (0.945,0.035) {$L$};
    \fi
    \ifnum\a>3
      \draw[hole] (1,0) rectangle (2,1);
      \node[text=wallred] at (1.5,0.5) {$3$};
    \fi
    \ifnum\a>6
      \draw[cell] (2,0) rectangle (3,1);
      \node at (2.5,0.5) {$6$};
      \node[anchor=south east,
            font=\fontsize{6}{6.5}\selectfont]
        at (2.945,0.035) {$R$};
    \fi
    \ifnum\a>8
      \draw[cell] (3,0) rectangle (4,1);
      \node at (3.5,0.5) {$1$};
      \node[anchor=south east,
            font=\fontsize{6}{6.5}\selectfont]
        at (3.945,0.035) {$L$};
    \fi
    \ifnum\a>9
      \draw[hole] (4,0) rectangle (5,1);
      \node[text=wallred] at (4.5,0.5) {$1$};
    \fi

    % Walls
    \ifnum\a>3
      \draw[wall] (1,0.05)--(1,0.95);
    \fi
    \ifnum\a>6
      \draw[wall] (2,0.05)--(2,0.95);
    \fi
    \ifnum\a>8
      \draw[wall] (3,0.05)--(3,0.95);
    \fi
    \ifnum\a>9
      \draw[wall] (4,0.05)--(4,0.95);
    \fi

    % Empty tableaux; omit the auxiliary star and its red label.
    \ifnum\a<3
      \node at (0,1.80) {$\varnothing$};
    \fi

  \end{scope}
  \end{scope}

  % Recorded pairs not yet used
  \def\pairs{\varnothing}

  \ifnum\a=1
    \def\pairs{\{(1,9)\}}
  \fi
  \ifnum\a=2
    \def\pairs{\{(1,9),(2,3)\}}
  \fi
  \ifnum\a=3
    \def\pairs{\{(1,9)\}}
  \fi
  \ifnum\a=4
    \def\pairs{\{(1,9)\}}
  \fi
  \ifnum\a=5
    \def\pairs{\{(1,9)\}}
  \fi
  \ifnum\a=6
    \def\pairs{\{(1,9),(6,7)\}}
  \fi
  \ifnum\a=7
    \def\pairs{\{(1,9)\}}
  \fi
  \ifnum\a=8
    \def\pairs{\{(1,9)\}}
  \fi

  \node[recorded] at (\panelx,{\panely-0.5})
    {$\mathcal R_{\a}=\pairs$};
}

% MOD 4: Arrow shaft length = twice this value.
\def\InverseArrowHalfLength{0.40} % total length = 0.80 cm
\foreach \a/\prev/\typ in {
  1/0/B, 2/1/B, 3/2/C,
  5/4/D, 6/5/B, 7/6/C,
  9/8/C, 10/9/A
}{
  \coordinate (stepcenter) at ($(stageright-\prev)!0.5!(stageleft-\a)$);
  \draw[step arrow]
    ($(stepcenter)+(-\InverseArrowHalfLength,0)$)--
    ($(stepcenter)+(\InverseArrowHalfLength,0)$)
    node[midway,above=2.4pt,font=\scriptsize]
      {$\textcolor{edgeblue}{e_{\a}}:\mathsf{\typ}$};
}
% MOD 5: 0.83 = 0.40 (half length) + 0.43 (gap to the tableau).

\foreach \a/\typ in {4/A,8/D}{
  \coordinate (stepcenter) at ($(stageleft-\a)+(-0.83,0)$);
  \draw[step arrow]
    ($(stepcenter)+(-\InverseArrowHalfLength,0)$)--
    ($(stepcenter)+(\InverseArrowHalfLength,0)$)
    node[midway,above=2.4pt,font=\scriptsize]
      {$\textcolor{edgeblue}{e_{\a}}:\mathsf{\typ}$};
}
\end{scope}
\end{tikzpicture}%
}
\vskip0.1in
\caption{An example of $\Psi^{-1}$, where  $\mathcal R_a$
lists the pairs recorded but not yet used after processing $e_a$. The initialization at $e_0$ gives $\rho(\star)=2$. \label{fig:Psi-inv}}
\end{figure}

\begin{enumerate}[
    label=\textup{(\arabic*)},
    leftmargin=2.5em,
    itemsep=4pt
]

\item \textbf{Initialization.}
Write \(e_0=(v,u)\). Since \(u\) is a leaf, set
\(\rho(\star)=\rho_{\mathrm \ell}(u)\).

\item 
For \(1\le a\le M\), assume that the entries corresponding to
\(e_1,\ldots,e_{a-1}\) have been processed. Inspect \(e_a\) and proceed
according to one of the following four cases.

\begin{enumerate}[
    leftmargin=3em,
    itemsep=3pt,
    topsep=3pt
]

\item
\uline{Suppose that \(e_a=(v,u)\) is of type \(\mathsf{A}\).}
Put \(a\) in the leftmost empty middle cell. If this cell is in
column \(j\), declare the bottom position in column \(j\) to be the
hole \(h_j\), and set
\(\rho(h_j)=\rho_{\mathrm \ell}(u)\).

\item
\uline{Suppose that $e_a=(v,w)$ is of type $\mathsf{B}$.}
Let $e_b$, with $b>a$, be the other incoming edge of $w$.
Record the pair $(a,b)$.

\item
\uline{Suppose that $e_a=(v,w)$ is of type $\mathsf{C}$.}
Let $(b,a)$ be the previously recorded pair.
Put $a$ in the leftmost empty middle cell and $b$ in the
bottom cell directly below it. Set
$\varepsilon(j)=\tau_w(e_a)$, where $j$ is the column containing $a$.

\item
\uline{Suppose that \(e_a\) is of type \(\mathsf{D}\).}
Put \(a\) in the leftmost empty top cell.

\end{enumerate}

\item
After \(e_M\) has been processed, the resulting decorated Young tableau with walls and holes is
defined to be \(\Psi^{-1}(N)\). 
\end{enumerate}
See Fig.~\ref{fig:Psi-inv} for an example of $\Psi^{-1}$. 

\begin{lemma}
\label{prop:inverse-well-defined}
The  map 
\(\Psi^{-1}:\mathcal N_{n,k,\ell} \to\widehat{\mathcal D}_{n,k,\ell}\)
is well defined.
\end{lemma}

\begin{proof}
Let $\Psi^{-1}(N)=(T;\varepsilon,\rho)$. Since $N$ has $k$ tree vertices and $\ell$ reticulation vertices, $T$ has $k$ top entries, $\ell$ bottom entries, and $n$ middle entries according to the construction of $\Psi^{-1}$. 
It is plain to see from the construction that the middle and top rows of $T$ are increasing, and every
occupied bottom entry is smaller than the middle entry above it.
Thus it remains only to verify the middle-to-top increasing  condition.

Let $m_1<m_2<\cdots<m_n$ and  \(t_1<\cdots<t_k\) be the middle and top
entries of $T$, respectively. Suppose that there are \(q_i\)
middle entries smaller than \(t_i\) in $T$. We need to verify that $q_i\geq i$ for any fixed $1\leq i\leq k$. Let 
\(E=(e_0,e_1,\ldots,e_M)\) be the canonical edge postorder of \(N\). 
According to the construction of $\Psi^{-1}$, there are exactly $q_i+1$ edges of types $\mathsf A$ or $\mathsf C$, and $i-1$ edges of type $\mathsf D$ in the prefix \(\pref_i(E):=(e_0,\ldots,e_{t_i-1})\) of $E$.
Suppose that there are  \(b_i\) edges of type $\mathsf B$ in the prefix $\pref_i(E)$. Thus, the prefix
\(\pref_i(E)\) contains exactly  \(q_i+b_i+i\) edges.

In edge postorder, every outgoing edge of a vertex precedes every
incoming edge of that vertex. Hence the prefix $\pref_i(E)$ contains both
outgoing edges of each of the \(i\) tree vertices entered by
\(e_{t_1},\ldots,e_{t_i}\), together with the unique outgoing
edge of each of the \(b_i\) reticulation vertices entered by
the type \(B\) edges. These are \(2i+b_i\) distinct
edges, since the contributing vertices are distinct. Therefore
\[
q_i+b_i+i\ge 2i+b_i,
\]
and hence \(q_i\ge i\). Thus at least \(i\) middle entries are
smaller than \(t_i\), giving \(m_i<t_i\), as required. This proves that $T$ is a Young tableau with walls and
holes in $\mathcal D_{n,k,\ell}$, completing the proof.  
\end{proof}

We are now in a position to prove Theorem~\ref{thm:tableau-network-correspondence}. 

\begin{proof}[{\bf Proof of Theorem~\ref{thm:tableau-network-correspondence}}]
Let \(X=(T;\varepsilon,\rho)\in
\widehat{\mathcal D}_{n,k,\ell}\). By Lemma~\ref{lem:edge-types}, the inverse construction in $\Psi^{-1}$ recovers
exactly the case in which each entry was processed in the 
construction of $\Psi$. The pairing of the two incoming edges at each
reticulation vertex  recovers the occupied bottom cells, while their local
labels recover \(\varepsilon\); similarly, \(e_0\) and the type \(\mathsf A\)
edges recover \(\rho\). Hence
\(\Psi^{-1}(\Psi(X))=X\) .

Conversely, starting from \(N\in\mathcal N_{n,k,\ell}\), the inverse
construction records from its canonical edge postorder exactly the
data used in the four local modifications defining
\(\Psi\). Applying the construction in $\Psi$ therefore
recovers the same edges, incidences, and local labels, in the same
canonical postorder. Thus
\(\Psi(\Psi^{-1}(N))=N\).

Therefore, \(\Psi\) and \(\Psi^{-1}\) are mutually inverse, and hence \(\Psi\) is a bijection.
\end{proof}

\medskip
\noindent\textbf{Declaration of AI Use.}
Generative AI was used to assist in finding our main involution $\Psi$ (the underlying rough ideas belong to the authors) and in drawing figures.

\section*{Acknowledgments}
We thank  Michael Fuchs for encouraging us to pursue a combinatorial proof of~\eqref{eq:intro-fuchs-yu} and Dongsu Kim for his helpful discussions on  Theorem~\ref{thm:weighted-symmetry} during their respective visits to Shandong University.
Z. Lin and J. Liu were  supported by the National Science Foundation of China (Grants 12322115 \& 12271301) and the Fundamental Research Funds for the Central Universities.
F.~Liu was supported by the Postdoctoral Fellowship Program and  China Postdoctoral Science Foundation (Grant No. BX2026002).

%\newpage

\clearpage
\newgeometry{top=0.65in,bottom=0.70in,left=0.70in,right=0.70in}
\appendix
\titleformat{\section}
  {\normalfont\Large\bfseries}
  {Appendix~\thesection.}{0.65em}{}
  
\section{Decorated binary networks of type \texorpdfstring{$(2,1,1)$}{(2,1,1)}}
\label{app:networks-211}
\vspace{-0.35em}

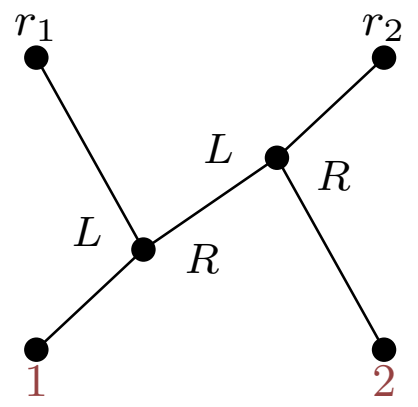
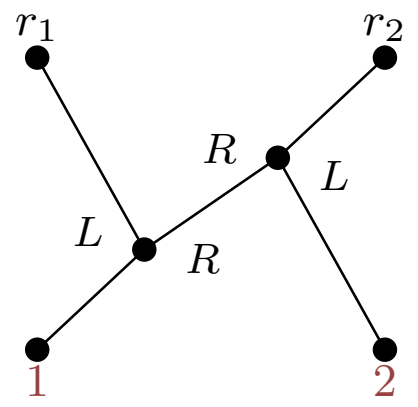
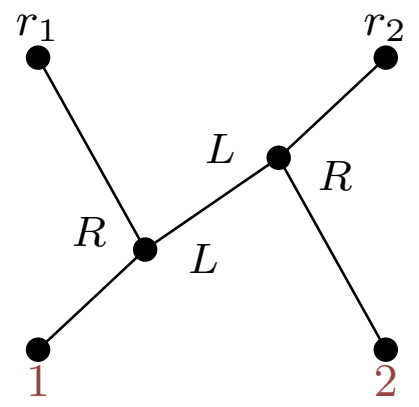
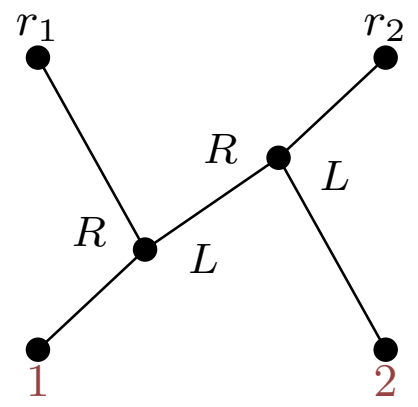
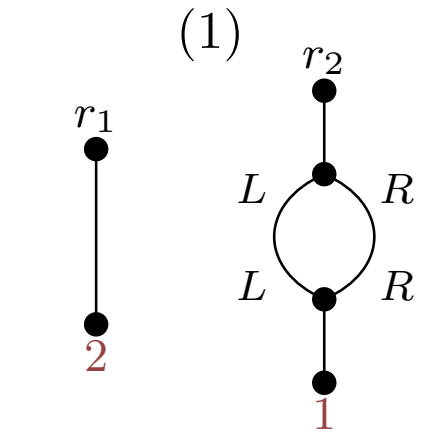
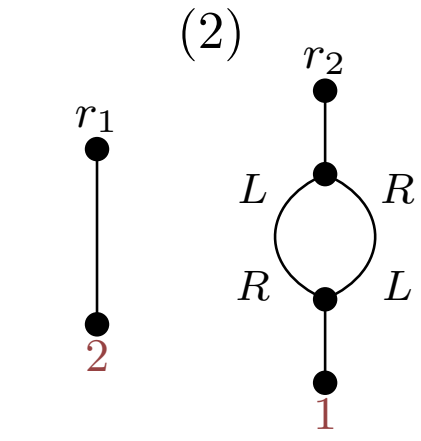
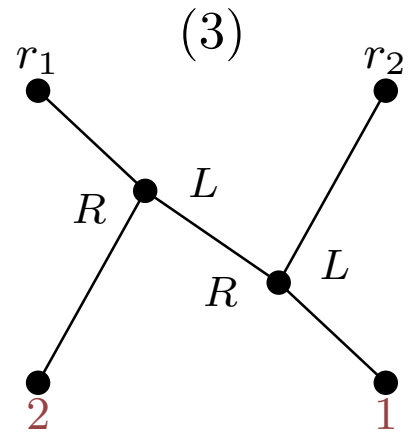
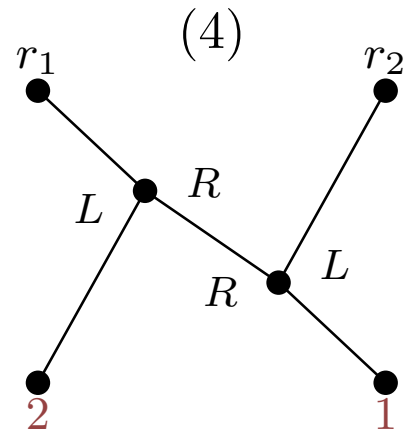
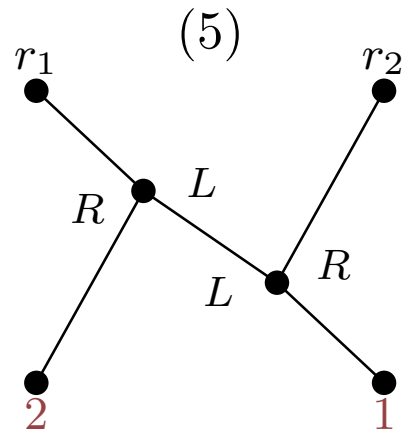
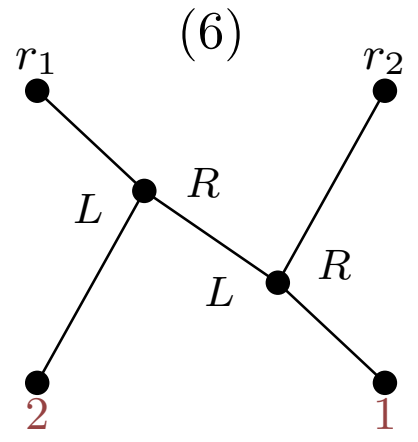
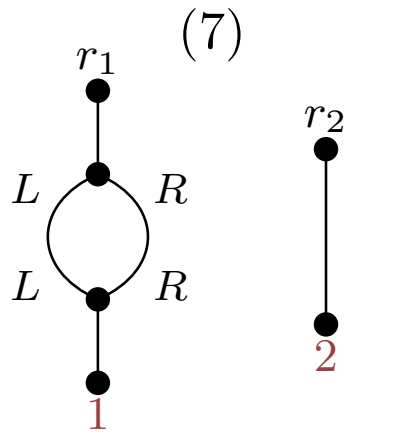
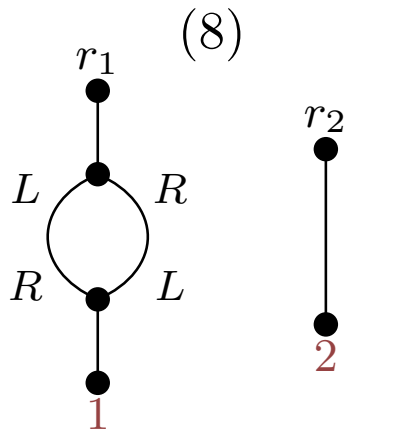
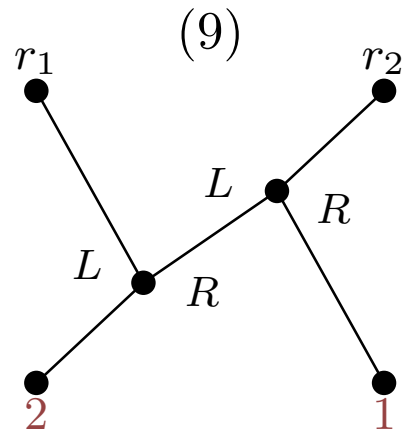
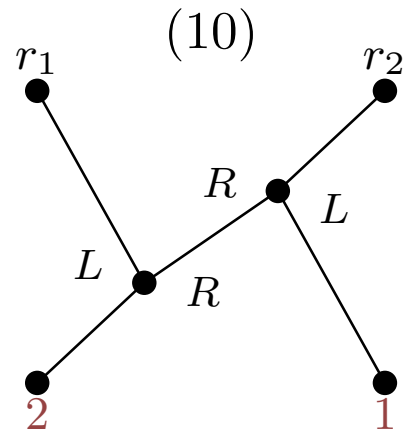
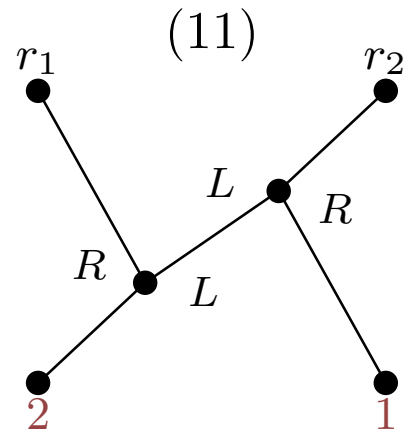
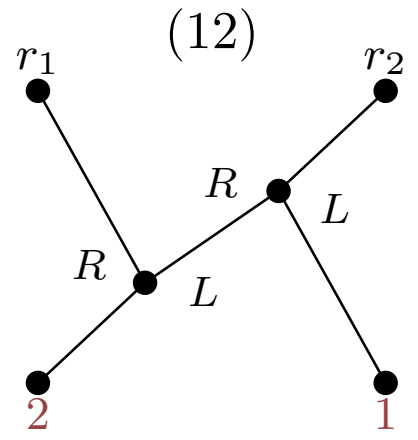
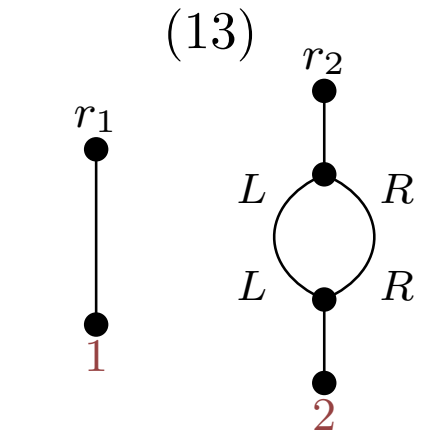
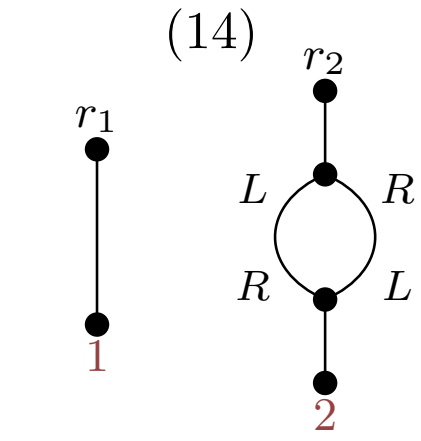
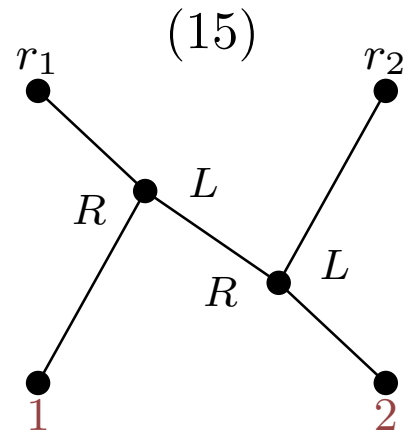
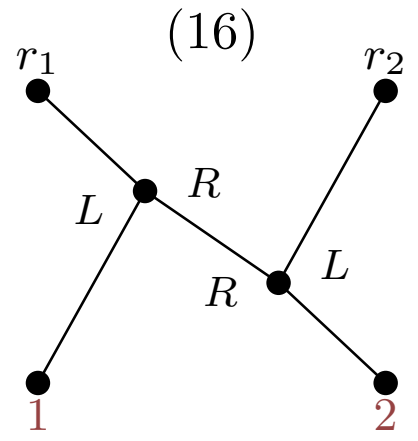
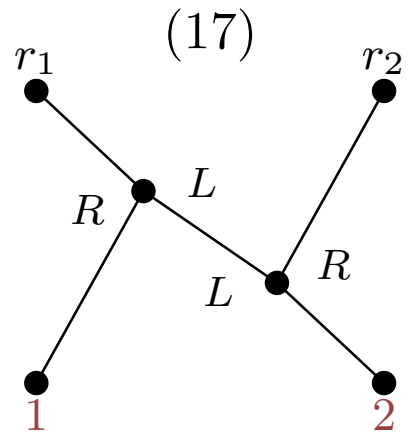
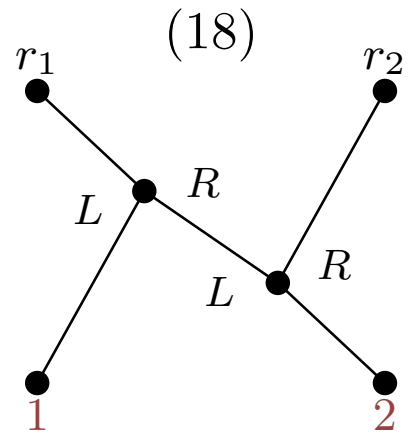
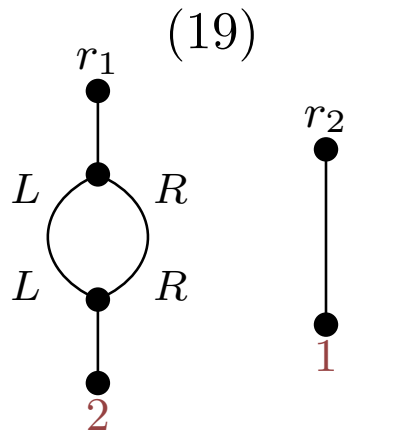
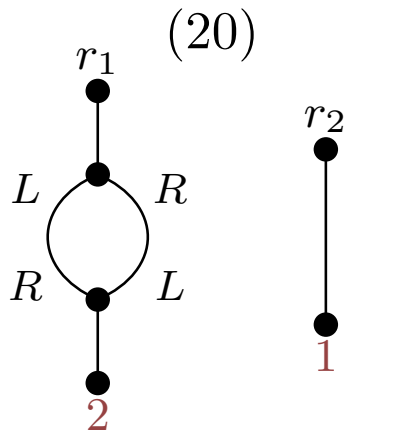
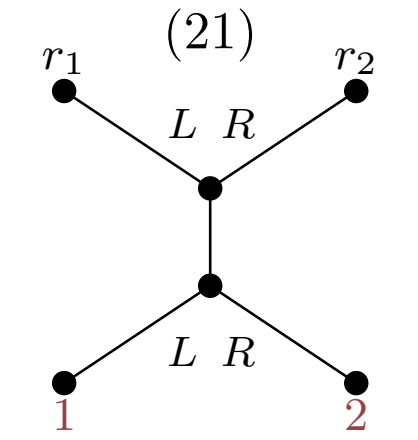
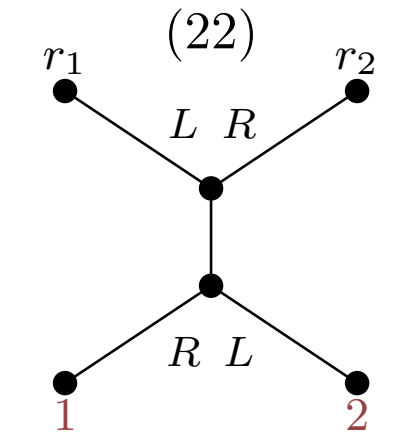
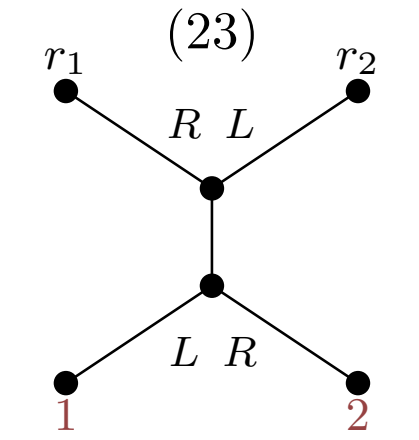
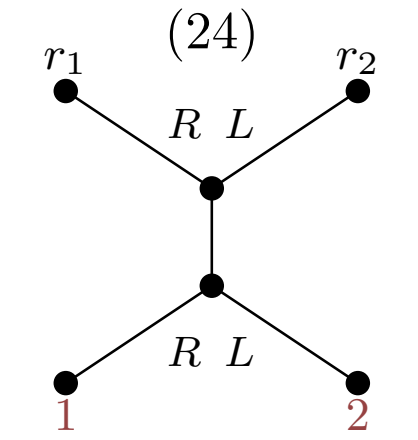
\tikzset{
  dbn edge/.style={draw=black,line width=0.49pt,line cap=round},
  dbn port/.style={font=\fontsize{8.0}{9.0}\selectfont,inner sep=0pt,outer sep=0pt},
  dbn boundary/.style={font=\fontsize{8.8}{9.8}\selectfont,inner sep=0pt,outer sep=0pt}
}

\newcommand{\dbncaption}[1]{\par\vspace{-0.15mm}{\footnotesize(#1)}}
\setlength{\tabcolsep}{0pt}
\noindent\makebox[\textwidth][c]{%
\begin{minipage}[t][3.18cm][t]{0.235\textwidth}\centering
% Class 1: same numbering, vertices, arcs and port orders as the uploaded HTML.
\begin{tikzpicture}[x=0.92cm,y=0.92cm]
  \path[use as bounding box] (-1.60,-0.55) rectangle (1.60,2.60);
  \coordinate (r1) at (-1.25000,2.10000);
  \coordinate (r2) at (1.25000,2.10000);
  \coordinate (t1) at (0.48000,1.38000);
  \coordinate (h1) at (-0.48000,0.72000);
  \coordinate (q1) at (-1.25000,0.00000);
  \coordinate (q2) at (1.25000,0.00000);
  % Arc 1: r1[-] -> h1[L].
  \draw[dbn edge] (r1) -- (h1);
  % Arc 2: r2[-] -> t1[-].
  \draw[dbn edge] (r2) -- (t1);
  % Arc 3: t1[R] -> q2[-].
  \draw[dbn edge] (t1) -- (q2);
  % Arc 4: t1[L] -> h1[R].
  \draw[dbn edge] (t1) -- (h1);
  % Arc 5: h1[-] -> q1[-].
  \draw[dbn edge] (h1) -- (q1);
  \fill[black] (r1) circle[radius=2.30pt];
  \fill[black] (r2) circle[radius=2.30pt];
  \fill[black] (t1) circle[radius=2.30pt];
  \fill[black] (h1) circle[radius=2.30pt];
  \fill[black] (q1) circle[radius=2.30pt];
  \fill[black] (q2) circle[radius=2.30pt];
  \node[dbn boundary,above=3.2pt] at (r1) {$r_1$};
  \node[dbn boundary,above=3.2pt] at (r2) {$r_2$};
  \node[dbn boundary,text=wallred,below=3.2pt] at (q1) {$1$};
  \node[dbn boundary,text=wallred,below=3.2pt] at (q2) {$2$};
  % Local in-port of h1 on arc 1.
  \node[dbn port] (p1-e1-in) at (-0.88429,0.84941) {$L$};
  % Local out-port of t1 on arc 3.
  \node[dbn port] (p1-e3-out) at (0.88429,1.25059) {$R$};
  % Local out-port of t1 on arc 4.
  \node[dbn port] (p1-e4-out) at (0.06025,1.44335) {$L$};
  % Local in-port of h1 on arc 4.
  \node[dbn port] (p1-e4-in) at (-0.06025,0.65665) {$R$};
\end{tikzpicture}
\dbncaption{1}
\end{minipage}%
\hspace{1.4mm}%
\begin{minipage}[t][3.18cm][t]{0.235\textwidth}\centering
% Class 2: same numbering, vertices, arcs and port orders as the uploaded HTML.
\begin{tikzpicture}[x=0.92cm,y=0.92cm]
  \path[use as bounding box] (-1.60,-0.55) rectangle (1.60,2.60);
  \coordinate (r1) at (-1.25000,2.10000);
  \coordinate (r2) at (1.25000,2.10000);
  \coordinate (t1) at (0.48000,1.38000);
  \coordinate (h1) at (-0.48000,0.72000);
  \coordinate (q1) at (-1.25000,0.00000);
  \coordinate (q2) at (1.25000,0.00000);
  % Arc 1: r1[-] -> h1[L].
  \draw[dbn edge] (r1) -- (h1);
  % Arc 2: r2[-] -> t1[-].
  \draw[dbn edge] (r2) -- (t1);
  % Arc 3: t1[L] -> q2[-].
  \draw[dbn edge] (t1) -- (q2);
  % Arc 4: t1[R] -> h1[R].
  \draw[dbn edge] (t1) -- (h1);
  % Arc 5: h1[-] -> q1[-].
  \draw[dbn edge] (h1) -- (q1);
  \fill[black] (r1) circle[radius=2.30pt];
  \fill[black] (r2) circle[radius=2.30pt];
  \fill[black] (t1) circle[radius=2.30pt];
  \fill[black] (h1) circle[radius=2.30pt];
  \fill[black] (q1) circle[radius=2.30pt];
  \fill[black] (q2) circle[radius=2.30pt];
  \node[dbn boundary,above=3.2pt] at (r1) {$r_1$};
  \node[dbn boundary,above=3.2pt] at (r2) {$r_2$};
  \node[dbn boundary,text=wallred,below=3.2pt] at (q1) {$1$};
  \node[dbn boundary,text=wallred,below=3.2pt] at (q2) {$2$};
  % Local in-port of h1 on arc 1.
  \node[dbn port] (p2-e1-in) at (-0.88429,0.84941) {$L$};
  % Local out-port of t1 on arc 3.
  \node[dbn port] (p2-e3-out) at (0.88429,1.25059) {$L$};
  % Local out-port of t1 on arc 4.
  \node[dbn port] (p2-e4-out) at (0.06025,1.44335) {$R$};
  % Local in-port of h1 on arc 4.
  \node[dbn port] (p2-e4-in) at (-0.06025,0.65665) {$R$};
\end{tikzpicture}
\dbncaption{2}
\end{minipage}%
\hspace{1.4mm}%
\begin{minipage}[t][3.18cm][t]{0.235\textwidth}\centering
% Class 3: same numbering, vertices, arcs and port orders as the uploaded HTML.
\begin{tikzpicture}[x=0.92cm,y=0.92cm]
  \path[use as bounding box] (-1.60,-0.55) rectangle (1.60,2.60);
  \coordinate (r1) at (-1.25000,2.10000);
  \coordinate (r2) at (1.25000,2.10000);
  \coordinate (t1) at (0.48000,1.38000);
  \coordinate (h1) at (-0.48000,0.72000);
  \coordinate (q1) at (-1.25000,0.00000);
  \coordinate (q2) at (1.25000,0.00000);
  % Arc 1: r1[-] -> h1[R].
  \draw[dbn edge] (r1) -- (h1);
  % Arc 2: r2[-] -> t1[-].
  \draw[dbn edge] (r2) -- (t1);
  % Arc 3: t1[R] -> q2[-].
  \draw[dbn edge] (t1) -- (q2);
  % Arc 4: t1[L] -> h1[L].
  \draw[dbn edge] (t1) -- (h1);
  % Arc 5: h1[-] -> q1[-].
  \draw[dbn edge] (h1) -- (q1);
  \fill[black] (r1) circle[radius=2.30pt];
  \fill[black] (r2) circle[radius=2.30pt];
  \fill[black] (t1) circle[radius=2.30pt];
  \fill[black] (h1) circle[radius=2.30pt];
  \fill[black] (q1) circle[radius=2.30pt];
  \fill[black] (q2) circle[radius=2.30pt];
  \node[dbn boundary,above=3.2pt] at (r1) {$r_1$};
  \node[dbn boundary,above=3.2pt] at (r2) {$r_2$};
  \node[dbn boundary,text=wallred,below=3.2pt] at (q1) {$1$};
  \node[dbn boundary,text=wallred,below=3.2pt] at (q2) {$2$};
  % Local in-port of h1 on arc 1.
  \node[dbn port] (p3-e1-in) at (-0.88429,0.84941) {$R$};
  % Local out-port of t1 on arc 3.
  \node[dbn port] (p3-e3-out) at (0.88429,1.25059) {$R$};
  % Local out-port of t1 on arc 4.
  \node[dbn port] (p3-e4-out) at (0.06025,1.44335) {$L$};
  % Local in-port of h1 on arc 4.
  \node[dbn port] (p3-e4-in) at (-0.06025,0.65665) {$L$};
\end{tikzpicture}
\dbncaption{3}
\end{minipage}%
\hspace{1.4mm}%
\begin{minipage}[t][3.18cm][t]{0.235\textwidth}\centering
% Class 4: same numbering, vertices, arcs and port orders as the uploaded HTML.
\begin{tikzpicture}[x=0.92cm,y=0.92cm]
  \path[use as bounding box] (-1.60,-0.55) rectangle (1.60,2.60);
  \coordinate (r1) at (-1.25000,2.10000);
  \coordinate (r2) at (1.25000,2.10000);
  \coordinate (t1) at (0.48000,1.38000);
  \coordinate (h1) at (-0.48000,0.72000);
  \coordinate (q1) at (-1.25000,0.00000);
  \coordinate (q2) at (1.25000,0.00000);
  % Arc 1: r1[-] -> h1[R].
  \draw[dbn edge] (r1) -- (h1);
  % Arc 2: r2[-] -> t1[-].
  \draw[dbn edge] (r2) -- (t1);
  % Arc 3: t1[L] -> q2[-].
  \draw[dbn edge] (t1) -- (q2);
  % Arc 4: t1[R] -> h1[L].
  \draw[dbn edge] (t1) -- (h1);
  % Arc 5: h1[-] -> q1[-].
  \draw[dbn edge] (h1) -- (q1);
  \fill[black] (r1) circle[radius=2.30pt];
  \fill[black] (r2) circle[radius=2.30pt];
  \fill[black] (t1) circle[radius=2.30pt];
  \fill[black] (h1) circle[radius=2.30pt];
  \fill[black] (q1) circle[radius=2.30pt];
  \fill[black] (q2) circle[radius=2.30pt];
  \node[dbn boundary,above=3.2pt] at (r1) {$r_1$};
  \node[dbn boundary,above=3.2pt] at (r2) {$r_2$};
  \node[dbn boundary,text=wallred,below=3.2pt] at (q1) {$1$};
  \node[dbn boundary,text=wallred,below=3.2pt] at (q2) {$2$};
  % Local in-port of h1 on arc 1.
  \node[dbn port] (p4-e1-in) at (-0.88429,0.84941) {$R$};
  % Local out-port of t1 on arc 3.
  \node[dbn port] (p4-e3-out) at (0.88429,1.25059) {$L$};
  % Local out-port of t1 on arc 4.
  \node[dbn port] (p4-e4-out) at (0.06025,1.44335) {$R$};
  % Local in-port of h1 on arc 4.
  \node[dbn port] (p4-e4-in) at (-0.06025,0.65665) {$L$};
\end{tikzpicture}
\dbncaption{4}
\end{minipage}%
}\par

\noindent\makebox[\textwidth][c]{%
\begin{minipage}[t][3.18cm][t]{0.235\textwidth}\centering
% Class 5: same numbering, vertices, arcs and port orders as the uploaded HTML.
\begin{tikzpicture}[x=0.92cm,y=0.92cm]
  \path[use as bounding box] (-1.60,-0.55) rectangle (1.60,2.60);
  \coordinate (r1) at (-0.82000,1.68000);
  \coordinate (r2) at (0.82000,2.10000);
  \coordinate (t1) at (0.82000,1.50000);
  \coordinate (h1) at (0.82000,0.60000);
  \coordinate (q1) at (0.82000,0.00000);
  \coordinate (q2) at (-0.82000,0.42000);
  % Arc 1: r1[-] -> q2[-].
  \draw[dbn edge] (r1) -- (q2);
  % Arc 2: r2[-] -> t1[-].
  \draw[dbn edge] (r2) -- (t1);
  % Arc 3: t1[L] -> h1[L].
  \draw[dbn edge] (t1) .. controls (0.34000,1.29000) and (0.34000,0.81000) .. (h1);
  % Arc 4: t1[R] -> h1[R].
  \draw[dbn edge] (t1) .. controls (1.30000,1.29000) and (1.30000,0.81000) .. (h1);
  % Arc 5: h1[-] -> q1[-].
  \draw[dbn edge] (h1) -- (q1);
  \fill[black] (r1) circle[radius=2.30pt];
  \fill[black] (r2) circle[radius=2.30pt];
  \fill[black] (t1) circle[radius=2.30pt];
  \fill[black] (h1) circle[radius=2.30pt];
  \fill[black] (q1) circle[radius=2.30pt];
  \fill[black] (q2) circle[radius=2.30pt];
  \node[dbn boundary,above=3.2pt] at (r1) {$r_1$};
  \node[dbn boundary,above=3.2pt] at (r2) {$r_2$};
  \node[dbn boundary,text=wallred,below=3.2pt] at (q1) {$1$};
  \node[dbn boundary,text=wallred,below=3.2pt] at (q2) {$2$};
  % Local out-port of t1 on arc 3.
  \node[dbn port] (p5-e3-out) at (0.30000,1.40000) {$L$};
  % Local in-port of h1 on arc 3.
  \node[dbn port] (p5-e3-in) at (0.30000,0.70000) {$L$};
  % Local out-port of t1 on arc 4.
  \node[dbn port] (p5-e4-out) at (1.34000,1.40000) {$R$};
  % Local in-port of h1 on arc 4.
  \node[dbn port] (p5-e4-in) at (1.34000,0.70000) {$R$};
\end{tikzpicture}
\dbncaption{5}
\end{minipage}%
\hspace{1.4mm}%
\begin{minipage}[t][3.18cm][t]{0.235\textwidth}\centering
% Class 6: same numbering, vertices, arcs and port orders as the uploaded HTML.
\begin{tikzpicture}[x=0.92cm,y=0.92cm]
  \path[use as bounding box] (-1.60,-0.55) rectangle (1.60,2.60);
  \coordinate (r1) at (-0.82000,1.68000);
  \coordinate (r2) at (0.82000,2.10000);
  \coordinate (t1) at (0.82000,1.50000);
  \coordinate (h1) at (0.82000,0.60000);
  \coordinate (q1) at (0.82000,0.00000);
  \coordinate (q2) at (-0.82000,0.42000);
  % Arc 1: r1[-] -> q2[-].
  \draw[dbn edge] (r1) -- (q2);
  % Arc 2: r2[-] -> t1[-].
  \draw[dbn edge] (r2) -- (t1);
  % Arc 3: t1[L] -> h1[R].
  \draw[dbn edge] (t1) .. controls (0.34000,1.29000) and (0.34000,0.81000) .. (h1);
  % Arc 4: t1[R] -> h1[L].
  \draw[dbn edge] (t1) .. controls (1.30000,1.29000) and (1.30000,0.81000) .. (h1);
  % Arc 5: h1[-] -> q1[-].
  \draw[dbn edge] (h1) -- (q1);
  \fill[black] (r1) circle[radius=2.30pt];
  \fill[black] (r2) circle[radius=2.30pt];
  \fill[black] (t1) circle[radius=2.30pt];
  \fill[black] (h1) circle[radius=2.30pt];
  \fill[black] (q1) circle[radius=2.30pt];
  \fill[black] (q2) circle[radius=2.30pt];
  \node[dbn boundary,above=3.2pt] at (r1) {$r_1$};
  \node[dbn boundary,above=3.2pt] at (r2) {$r_2$};
  \node[dbn boundary,text=wallred,below=3.2pt] at (q1) {$1$};
  \node[dbn boundary,text=wallred,below=3.2pt] at (q2) {$2$};
  % Local out-port of t1 on arc 3.
  \node[dbn port] (p6-e3-out) at (0.30000,1.40000) {$L$};
  % Local in-port of h1 on arc 3.
  \node[dbn port] (p6-e3-in) at (0.30000,0.70000) {$R$};
  % Local out-port of t1 on arc 4.
  \node[dbn port] (p6-e4-out) at (1.34000,1.40000) {$R$};
  % Local in-port of h1 on arc 4.
  \node[dbn port] (p6-e4-in) at (1.34000,0.70000) {$L$};
\end{tikzpicture}
\dbncaption{6}
\end{minipage}%
\hspace{1.4mm}%
\begin{minipage}[t][3.18cm][t]{0.235\textwidth}\centering
% Class 7: same numbering, vertices, arcs and port orders as the uploaded HTML.
\begin{tikzpicture}[x=0.92cm,y=0.92cm]
  \path[use as bounding box] (-1.60,-0.55) rectangle (1.60,2.60);
  \coordinate (r1) at (-1.25000,2.10000);
  \coordinate (r2) at (1.25000,2.10000);
  \coordinate (t1) at (-0.48000,1.38000);
  \coordinate (h1) at (0.48000,0.72000);
  \coordinate (q1) at (1.25000,0.00000);
  \coordinate (q2) at (-1.25000,0.00000);
  % Arc 1: r1[-] -> t1[-].
  \draw[dbn edge] (r1) -- (t1);
  % Arc 2: r2[-] -> h1[L].
  \draw[dbn edge] (r2) -- (h1);
  % Arc 3: t1[R] -> q2[-].
  \draw[dbn edge] (t1) -- (q2);
  % Arc 4: t1[L] -> h1[R].
  \draw[dbn edge] (t1) -- (h1);
  % Arc 5: h1[-] -> q1[-].
  \draw[dbn edge] (h1) -- (q1);
  \fill[black] (r1) circle[radius=2.30pt];
  \fill[black] (r2) circle[radius=2.30pt];
  \fill[black] (t1) circle[radius=2.30pt];
  \fill[black] (h1) circle[radius=2.30pt];
  \fill[black] (q1) circle[radius=2.30pt];
  \fill[black] (q2) circle[radius=2.30pt];
  \node[dbn boundary,above=3.2pt] at (r1) {$r_1$};
  \node[dbn boundary,above=3.2pt] at (r2) {$r_2$};
  \node[dbn boundary,text=wallred,below=3.2pt] at (q1) {$1$};
  \node[dbn boundary,text=wallred,below=3.2pt] at (q2) {$2$};
  % Local in-port of h1 on arc 2.
  \node[dbn port] (p7-e2-in) at (0.88429,0.84941) {$L$};
  % Local out-port of t1 on arc 3.
  \node[dbn port] (p7-e3-out) at (-0.88429,1.25059) {$R$};
  % Local out-port of t1 on arc 4.
  \node[dbn port] (p7-e4-out) at (-0.06025,1.44335) {$L$};
  % Local in-port of h1 on arc 4.
  \node[dbn port] (p7-e4-in) at (0.06025,0.65665) {$R$};
\end{tikzpicture}
\dbncaption{7}
\end{minipage}%
\hspace{1.4mm}%
\begin{minipage}[t][3.18cm][t]{0.235\textwidth}\centering
% Class 8: same numbering, vertices, arcs and port orders as the uploaded HTML.
\begin{tikzpicture}[x=0.92cm,y=0.92cm]
  \path[use as bounding box] (-1.60,-0.55) rectangle (1.60,2.60);
  \coordinate (r1) at (-1.25000,2.10000);
  \coordinate (r2) at (1.25000,2.10000);
  \coordinate (t1) at (-0.48000,1.38000);
  \coordinate (h1) at (0.48000,0.72000);
  \coordinate (q1) at (1.25000,0.00000);
  \coordinate (q2) at (-1.25000,0.00000);
  % Arc 1: r1[-] -> t1[-].
  \draw[dbn edge] (r1) -- (t1);
  % Arc 2: r2[-] -> h1[L].
  \draw[dbn edge] (r2) -- (h1);
  % Arc 3: t1[L] -> q2[-].
  \draw[dbn edge] (t1) -- (q2);
  % Arc 4: t1[R] -> h1[R].
  \draw[dbn edge] (t1) -- (h1);
  % Arc 5: h1[-] -> q1[-].
  \draw[dbn edge] (h1) -- (q1);
  \fill[black] (r1) circle[radius=2.30pt];
  \fill[black] (r2) circle[radius=2.30pt];
  \fill[black] (t1) circle[radius=2.30pt];
  \fill[black] (h1) circle[radius=2.30pt];
  \fill[black] (q1) circle[radius=2.30pt];
  \fill[black] (q2) circle[radius=2.30pt];
  \node[dbn boundary,above=3.2pt] at (r1) {$r_1$};
  \node[dbn boundary,above=3.2pt] at (r2) {$r_2$};
  \node[dbn boundary,text=wallred,below=3.2pt] at (q1) {$1$};
  \node[dbn boundary,text=wallred,below=3.2pt] at (q2) {$2$};
  % Local in-port of h1 on arc 2.
  \node[dbn port] (p8-e2-in) at (0.88429,0.84941) {$L$};
  % Local out-port of t1 on arc 3.
  \node[dbn port] (p8-e3-out) at (-0.88429,1.25059) {$L$};
  % Local out-port of t1 on arc 4.
  \node[dbn port] (p8-e4-out) at (-0.06025,1.44335) {$R$};
  % Local in-port of h1 on arc 4.
  \node[dbn port] (p8-e4-in) at (0.06025,0.65665) {$R$};
\end{tikzpicture}
\dbncaption{8}
\end{minipage}%
}\par

\noindent\makebox[\textwidth][c]{%
\begin{minipage}[t][3.18cm][t]{0.235\textwidth}\centering
% Class 9: same numbering, vertices, arcs and port orders as the uploaded HTML.
\begin{tikzpicture}[x=0.92cm,y=0.92cm]
  \path[use as bounding box] (-1.60,-0.55) rectangle (1.60,2.60);
  \coordinate (r1) at (-1.25000,2.10000);
  \coordinate (r2) at (1.25000,2.10000);
  \coordinate (t1) at (-0.48000,1.38000);
  \coordinate (h1) at (0.48000,0.72000);
  \coordinate (q1) at (1.25000,0.00000);
  \coordinate (q2) at (-1.25000,0.00000);
  % Arc 1: r1[-] -> t1[-].
  \draw[dbn edge] (r1) -- (t1);
  % Arc 2: r2[-] -> h1[R].
  \draw[dbn edge] (r2) -- (h1);
  % Arc 3: t1[R] -> q2[-].
  \draw[dbn edge] (t1) -- (q2);
  % Arc 4: t1[L] -> h1[L].
  \draw[dbn edge] (t1) -- (h1);
  % Arc 5: h1[-] -> q1[-].
  \draw[dbn edge] (h1) -- (q1);
  \fill[black] (r1) circle[radius=2.30pt];
  \fill[black] (r2) circle[radius=2.30pt];
  \fill[black] (t1) circle[radius=2.30pt];
  \fill[black] (h1) circle[radius=2.30pt];
  \fill[black] (q1) circle[radius=2.30pt];
  \fill[black] (q2) circle[radius=2.30pt];
  \node[dbn boundary,above=3.2pt] at (r1) {$r_1$};
  \node[dbn boundary,above=3.2pt] at (r2) {$r_2$};
  \node[dbn boundary,text=wallred,below=3.2pt] at (q1) {$1$};
  \node[dbn boundary,text=wallred,below=3.2pt] at (q2) {$2$};
  % Local in-port of h1 on arc 2.
  \node[dbn port] (p9-e2-in) at (0.88429,0.84941) {$R$};
  % Local out-port of t1 on arc 3.
  \node[dbn port] (p9-e3-out) at (-0.88429,1.25059) {$R$};
  % Local out-port of t1 on arc 4.
  \node[dbn port] (p9-e4-out) at (-0.06025,1.44335) {$L$};
  % Local in-port of h1 on arc 4.
  \node[dbn port] (p9-e4-in) at (0.06025,0.65665) {$L$};
\end{tikzpicture}
\dbncaption{9}
\end{minipage}%
\hspace{1.4mm}%
\begin{minipage}[t][3.18cm][t]{0.235\textwidth}\centering
% Class 10: same numbering, vertices, arcs and port orders as the uploaded HTML.
\begin{tikzpicture}[x=0.92cm,y=0.92cm]
  \path[use as bounding box] (-1.60,-0.55) rectangle (1.60,2.60);
  \coordinate (r1) at (-1.25000,2.10000);
  \coordinate (r2) at (1.25000,2.10000);
  \coordinate (t1) at (-0.48000,1.38000);
  \coordinate (h1) at (0.48000,0.72000);
  \coordinate (q1) at (1.25000,0.00000);
  \coordinate (q2) at (-1.25000,0.00000);
  % Arc 1: r1[-] -> t1[-].
  \draw[dbn edge] (r1) -- (t1);
  % Arc 2: r2[-] -> h1[R].
  \draw[dbn edge] (r2) -- (h1);
  % Arc 3: t1[L] -> q2[-].
  \draw[dbn edge] (t1) -- (q2);
  % Arc 4: t1[R] -> h1[L].
  \draw[dbn edge] (t1) -- (h1);
  % Arc 5: h1[-] -> q1[-].
  \draw[dbn edge] (h1) -- (q1);
  \fill[black] (r1) circle[radius=2.30pt];
  \fill[black] (r2) circle[radius=2.30pt];
  \fill[black] (t1) circle[radius=2.30pt];
  \fill[black] (h1) circle[radius=2.30pt];
  \fill[black] (q1) circle[radius=2.30pt];
  \fill[black] (q2) circle[radius=2.30pt];
  \node[dbn boundary,above=3.2pt] at (r1) {$r_1$};
  \node[dbn boundary,above=3.2pt] at (r2) {$r_2$};
  \node[dbn boundary,text=wallred,below=3.2pt] at (q1) {$1$};
  \node[dbn boundary,text=wallred,below=3.2pt] at (q2) {$2$};
  % Local in-port of h1 on arc 2.
  \node[dbn port] (p10-e2-in) at (0.88429,0.84941) {$R$};
  % Local out-port of t1 on arc 3.
  \node[dbn port] (p10-e3-out) at (-0.88429,1.25059) {$L$};
  % Local out-port of t1 on arc 4.
  \node[dbn port] (p10-e4-out) at (-0.06025,1.44335) {$R$};
  % Local in-port of h1 on arc 4.
  \node[dbn port] (p10-e4-in) at (0.06025,0.65665) {$L$};
\end{tikzpicture}
\dbncaption{10}
\end{minipage}%
\hspace{1.4mm}%
\begin{minipage}[t][3.18cm][t]{0.235\textwidth}\centering
% Class 11: same numbering, vertices, arcs and port orders as the uploaded HTML.
\begin{tikzpicture}[x=0.92cm,y=0.92cm]
  \path[use as bounding box] (-1.60,-0.55) rectangle (1.60,2.60);
  \coordinate (r1) at (-0.82000,2.10000);
  \coordinate (r2) at (0.82000,1.68000);
  \coordinate (t1) at (-0.82000,1.50000);
  \coordinate (h1) at (-0.82000,0.60000);
  \coordinate (q1) at (-0.82000,0.00000);
  \coordinate (q2) at (0.82000,0.42000);
  % Arc 1: r1[-] -> t1[-].
  \draw[dbn edge] (r1) -- (t1);
  % Arc 2: r2[-] -> q2[-].
  \draw[dbn edge] (r2) -- (q2);
  % Arc 3: t1[L] -> h1[L].
  \draw[dbn edge] (t1) .. controls (-1.30000,1.29000) and (-1.30000,0.81000) .. (h1);
  % Arc 4: t1[R] -> h1[R].
  \draw[dbn edge] (t1) .. controls (-0.34000,1.29000) and (-0.34000,0.81000) .. (h1);
  % Arc 5: h1[-] -> q1[-].
  \draw[dbn edge] (h1) -- (q1);
  \fill[black] (r1) circle[radius=2.30pt];
  \fill[black] (r2) circle[radius=2.30pt];
  \fill[black] (t1) circle[radius=2.30pt];
  \fill[black] (h1) circle[radius=2.30pt];
  \fill[black] (q1) circle[radius=2.30pt];
  \fill[black] (q2) circle[radius=2.30pt];
  \node[dbn boundary,above=3.2pt] at (r1) {$r_1$};
  \node[dbn boundary,above=3.2pt] at (r2) {$r_2$};
  \node[dbn boundary,text=wallred,below=3.2pt] at (q1) {$1$};
  \node[dbn boundary,text=wallred,below=3.2pt] at (q2) {$2$};
  % Local out-port of t1 on arc 3.
  \node[dbn port] (p11-e3-out) at (-1.34000,1.40000) {$L$};
  % Local in-port of h1 on arc 3.
  \node[dbn port] (p11-e3-in) at (-1.34000,0.70000) {$L$};
  % Local out-port of t1 on arc 4.
  \node[dbn port] (p11-e4-out) at (-0.30000,1.40000) {$R$};
  % Local in-port of h1 on arc 4.
  \node[dbn port] (p11-e4-in) at (-0.30000,0.70000) {$R$};
\end{tikzpicture}
\dbncaption{11}
\end{minipage}%
\hspace{1.4mm}%
\begin{minipage}[t][3.18cm][t]{0.235\textwidth}\centering
% Class 12: same numbering, vertices, arcs and port orders as the uploaded HTML.
\begin{tikzpicture}[x=0.92cm,y=0.92cm]
  \path[use as bounding box] (-1.60,-0.55) rectangle (1.60,2.60);
  \coordinate (r1) at (-0.82000,2.10000);
  \coordinate (r2) at (0.82000,1.68000);
  \coordinate (t1) at (-0.82000,1.50000);
  \coordinate (h1) at (-0.82000,0.60000);
  \coordinate (q1) at (-0.82000,0.00000);
  \coordinate (q2) at (0.82000,0.42000);
  % Arc 1: r1[-] -> t1[-].
  \draw[dbn edge] (r1) -- (t1);
  % Arc 2: r2[-] -> q2[-].
  \draw[dbn edge] (r2) -- (q2);
  % Arc 3: t1[L] -> h1[R].
  \draw[dbn edge] (t1) .. controls (-1.30000,1.29000) and (-1.30000,0.81000) .. (h1);
  % Arc 4: t1[R] -> h1[L].
  \draw[dbn edge] (t1) .. controls (-0.34000,1.29000) and (-0.34000,0.81000) .. (h1);
  % Arc 5: h1[-] -> q1[-].
  \draw[dbn edge] (h1) -- (q1);
  \fill[black] (r1) circle[radius=2.30pt];
  \fill[black] (r2) circle[radius=2.30pt];
  \fill[black] (t1) circle[radius=2.30pt];
  \fill[black] (h1) circle[radius=2.30pt];
  \fill[black] (q1) circle[radius=2.30pt];
  \fill[black] (q2) circle[radius=2.30pt];
  \node[dbn boundary,above=3.2pt] at (r1) {$r_1$};
  \node[dbn boundary,above=3.2pt] at (r2) {$r_2$};
  \node[dbn boundary,text=wallred,below=3.2pt] at (q1) {$1$};
  \node[dbn boundary,text=wallred,below=3.2pt] at (q2) {$2$};
  % Local out-port of t1 on arc 3.
  \node[dbn port] (p12-e3-out) at (-1.34000,1.40000) {$L$};
  % Local in-port of h1 on arc 3.
  \node[dbn port] (p12-e3-in) at (-1.34000,0.70000) {$R$};
  % Local out-port of t1 on arc 4.
  \node[dbn port] (p12-e4-out) at (-0.30000,1.40000) {$R$};
  % Local in-port of h1 on arc 4.
  \node[dbn port] (p12-e4-in) at (-0.30000,0.70000) {$L$};
\end{tikzpicture}
\dbncaption{12}
\end{minipage}%
}\par

\noindent\makebox[\textwidth][c]{%
\begin{minipage}[t][3.18cm][t]{0.235\textwidth}\centering
% Class 13: same numbering, vertices, arcs and port orders as the uploaded HTML.
\begin{tikzpicture}[x=0.92cm,y=0.92cm]
  \path[use as bounding box] (-1.60,-0.55) rectangle (1.60,2.60);
  \coordinate (r1) at (-1.25000,2.10000);
  \coordinate (r2) at (1.25000,2.10000);
  \coordinate (t1) at (0.48000,1.38000);
  \coordinate (h1) at (-0.48000,0.72000);
  \coordinate (q1) at (1.25000,0.00000);
  \coordinate (q2) at (-1.25000,0.00000);
  % Arc 1: r1[-] -> h1[L].
  \draw[dbn edge] (r1) -- (h1);
  % Arc 2: r2[-] -> t1[-].
  \draw[dbn edge] (r2) -- (t1);
  % Arc 3: t1[R] -> q1[-].
  \draw[dbn edge] (t1) -- (q1);
  % Arc 4: t1[L] -> h1[R].
  \draw[dbn edge] (t1) -- (h1);
  % Arc 5: h1[-] -> q2[-].
  \draw[dbn edge] (h1) -- (q2);
  \fill[black] (r1) circle[radius=2.30pt];
  \fill[black] (r2) circle[radius=2.30pt];
  \fill[black] (t1) circle[radius=2.30pt];
  \fill[black] (h1) circle[radius=2.30pt];
  \fill[black] (q1) circle[radius=2.30pt];
  \fill[black] (q2) circle[radius=2.30pt];
  \node[dbn boundary,above=3.2pt] at (r1) {$r_1$};
  \node[dbn boundary,above=3.2pt] at (r2) {$r_2$};
  \node[dbn boundary,text=wallred,below=3.2pt] at (q1) {$1$};
  \node[dbn boundary,text=wallred,below=3.2pt] at (q2) {$2$};
  % Local in-port of h1 on arc 1.
  \node[dbn port] (p13-e1-in) at (-0.88429,0.84941) {$L$};
  % Local out-port of t1 on arc 3.
  \node[dbn port] (p13-e3-out) at (0.88429,1.25059) {$R$};
  % Local out-port of t1 on arc 4.
  \node[dbn port] (p13-e4-out) at (0.06025,1.44335) {$L$};
  % Local in-port of h1 on arc 4.
  \node[dbn port] (p13-e4-in) at (-0.06025,0.65665) {$R$};
\end{tikzpicture}
\dbncaption{13}
\end{minipage}%
\hspace{1.4mm}%
\begin{minipage}[t][3.18cm][t]{0.235\textwidth}\centering
% Class 14: same numbering, vertices, arcs and port orders as the uploaded HTML.
\begin{tikzpicture}[x=0.92cm,y=0.92cm]
  \path[use as bounding box] (-1.60,-0.55) rectangle (1.60,2.60);
  \coordinate (r1) at (-1.25000,2.10000);
  \coordinate (r2) at (1.25000,2.10000);
  \coordinate (t1) at (0.48000,1.38000);
  \coordinate (h1) at (-0.48000,0.72000);
  \coordinate (q1) at (1.25000,0.00000);
  \coordinate (q2) at (-1.25000,0.00000);
  % Arc 1: r1[-] -> h1[L].
  \draw[dbn edge] (r1) -- (h1);
  % Arc 2: r2[-] -> t1[-].
  \draw[dbn edge] (r2) -- (t1);
  % Arc 3: t1[L] -> q1[-].
  \draw[dbn edge] (t1) -- (q1);
  % Arc 4: t1[R] -> h1[R].
  \draw[dbn edge] (t1) -- (h1);
  % Arc 5: h1[-] -> q2[-].
  \draw[dbn edge] (h1) -- (q2);
  \fill[black] (r1) circle[radius=2.30pt];
  \fill[black] (r2) circle[radius=2.30pt];
  \fill[black] (t1) circle[radius=2.30pt];
  \fill[black] (h1) circle[radius=2.30pt];
  \fill[black] (q1) circle[radius=2.30pt];
  \fill[black] (q2) circle[radius=2.30pt];
  \node[dbn boundary,above=3.2pt] at (r1) {$r_1$};
  \node[dbn boundary,above=3.2pt] at (r2) {$r_2$};
  \node[dbn boundary,text=wallred,below=3.2pt] at (q1) {$1$};
  \node[dbn boundary,text=wallred,below=3.2pt] at (q2) {$2$};
  % Local in-port of h1 on arc 1.
  \node[dbn port] (p14-e1-in) at (-0.88429,0.84941) {$L$};
  % Local out-port of t1 on arc 3.
  \node[dbn port] (p14-e3-out) at (0.88429,1.25059) {$L$};
  % Local out-port of t1 on arc 4.
  \node[dbn port] (p14-e4-out) at (0.06025,1.44335) {$R$};
  % Local in-port of h1 on arc 4.
  \node[dbn port] (p14-e4-in) at (-0.06025,0.65665) {$R$};
\end{tikzpicture}
\dbncaption{14}
\end{minipage}%
\hspace{1.4mm}%
\begin{minipage}[t][3.18cm][t]{0.235\textwidth}\centering
% Class 15: same numbering, vertices, arcs and port orders as the uploaded HTML.
\begin{tikzpicture}[x=0.92cm,y=0.92cm]
  \path[use as bounding box] (-1.60,-0.55) rectangle (1.60,2.60);
  \coordinate (r1) at (-1.25000,2.10000);
  \coordinate (r2) at (1.25000,2.10000);
  \coordinate (t1) at (0.48000,1.38000);
  \coordinate (h1) at (-0.48000,0.72000);
  \coordinate (q1) at (1.25000,0.00000);
  \coordinate (q2) at (-1.25000,0.00000);
  % Arc 1: r1[-] -> h1[R].
  \draw[dbn edge] (r1) -- (h1);
  % Arc 2: r2[-] -> t1[-].
  \draw[dbn edge] (r2) -- (t1);
  % Arc 3: t1[R] -> q1[-].
  \draw[dbn edge] (t1) -- (q1);
  % Arc 4: t1[L] -> h1[L].
  \draw[dbn edge] (t1) -- (h1);
  % Arc 5: h1[-] -> q2[-].
  \draw[dbn edge] (h1) -- (q2);
  \fill[black] (r1) circle[radius=2.30pt];
  \fill[black] (r2) circle[radius=2.30pt];
  \fill[black] (t1) circle[radius=2.30pt];
  \fill[black] (h1) circle[radius=2.30pt];
  \fill[black] (q1) circle[radius=2.30pt];
  \fill[black] (q2) circle[radius=2.30pt];
  \node[dbn boundary,above=3.2pt] at (r1) {$r_1$};
  \node[dbn boundary,above=3.2pt] at (r2) {$r_2$};
  \node[dbn boundary,text=wallred,below=3.2pt] at (q1) {$1$};
  \node[dbn boundary,text=wallred,below=3.2pt] at (q2) {$2$};
  % Local in-port of h1 on arc 1.
  \node[dbn port] (p15-e1-in) at (-0.88429,0.84941) {$R$};
  % Local out-port of t1 on arc 3.
  \node[dbn port] (p15-e3-out) at (0.88429,1.25059) {$R$};
  % Local out-port of t1 on arc 4.
  \node[dbn port] (p15-e4-out) at (0.06025,1.44335) {$L$};
  % Local in-port of h1 on arc 4.
  \node[dbn port] (p15-e4-in) at (-0.06025,0.65665) {$L$};
\end{tikzpicture}
\dbncaption{15}
\end{minipage}%
\hspace{1.4mm}%
\begin{minipage}[t][3.18cm][t]{0.235\textwidth}\centering
% Class 16: same numbering, vertices, arcs and port orders as the uploaded HTML.
\begin{tikzpicture}[x=0.92cm,y=0.92cm]
  \path[use as bounding box] (-1.60,-0.55) rectangle (1.60,2.60);
  \coordinate (r1) at (-1.25000,2.10000);
  \coordinate (r2) at (1.25000,2.10000);
  \coordinate (t1) at (0.48000,1.38000);
  \coordinate (h1) at (-0.48000,0.72000);
  \coordinate (q1) at (1.25000,0.00000);
  \coordinate (q2) at (-1.25000,0.00000);
  % Arc 1: r1[-] -> h1[R].
  \draw[dbn edge] (r1) -- (h1);
  % Arc 2: r2[-] -> t1[-].
  \draw[dbn edge] (r2) -- (t1);
  % Arc 3: t1[L] -> q1[-].
  \draw[dbn edge] (t1) -- (q1);
  % Arc 4: t1[R] -> h1[L].
  \draw[dbn edge] (t1) -- (h1);
  % Arc 5: h1[-] -> q2[-].
  \draw[dbn edge] (h1) -- (q2);
  \fill[black] (r1) circle[radius=2.30pt];
  \fill[black] (r2) circle[radius=2.30pt];
  \fill[black] (t1) circle[radius=2.30pt];
  \fill[black] (h1) circle[radius=2.30pt];
  \fill[black] (q1) circle[radius=2.30pt];
  \fill[black] (q2) circle[radius=2.30pt];
  \node[dbn boundary,above=3.2pt] at (r1) {$r_1$};
  \node[dbn boundary,above=3.2pt] at (r2) {$r_2$};
  \node[dbn boundary,text=wallred,below=3.2pt] at (q1) {$1$};
  \node[dbn boundary,text=wallred,below=3.2pt] at (q2) {$2$};
  % Local in-port of h1 on arc 1.
  \node[dbn port] (p16-e1-in) at (-0.88429,0.84941) {$R$};
  % Local out-port of t1 on arc 3.
  \node[dbn port] (p16-e3-out) at (0.88429,1.25059) {$L$};
  % Local out-port of t1 on arc 4.
  \node[dbn port] (p16-e4-out) at (0.06025,1.44335) {$R$};
  % Local in-port of h1 on arc 4.
  \node[dbn port] (p16-e4-in) at (-0.06025,0.65665) {$L$};
\end{tikzpicture}
\dbncaption{16}
\end{minipage}%
}\par

\noindent\makebox[\textwidth][c]{%
\begin{minipage}[t][3.18cm][t]{0.235\textwidth}\centering
% Class 17: same numbering, vertices, arcs and port orders as the uploaded HTML.
\begin{tikzpicture}[x=0.92cm,y=0.92cm]
  \path[use as bounding box] (-1.60,-0.55) rectangle (1.60,2.60);
  \coordinate (r1) at (-0.82000,1.68000);
  \coordinate (r2) at (0.82000,2.10000);
  \coordinate (t1) at (0.82000,1.50000);
  \coordinate (h1) at (0.82000,0.60000);
  \coordinate (q1) at (-0.82000,0.42000);
  \coordinate (q2) at (0.82000,0.00000);
  % Arc 1: r1[-] -> q1[-].
  \draw[dbn edge] (r1) -- (q1);
  % Arc 2: r2[-] -> t1[-].
  \draw[dbn edge] (r2) -- (t1);
  % Arc 3: t1[L] -> h1[L].
  \draw[dbn edge] (t1) .. controls (0.34000,1.29000) and (0.34000,0.81000) .. (h1);
  % Arc 4: t1[R] -> h1[R].
  \draw[dbn edge] (t1) .. controls (1.30000,1.29000) and (1.30000,0.81000) .. (h1);
  % Arc 5: h1[-] -> q2[-].
  \draw[dbn edge] (h1) -- (q2);
  \fill[black] (r1) circle[radius=2.30pt];
  \fill[black] (r2) circle[radius=2.30pt];
  \fill[black] (t1) circle[radius=2.30pt];
  \fill[black] (h1) circle[radius=2.30pt];
  \fill[black] (q1) circle[radius=2.30pt];
  \fill[black] (q2) circle[radius=2.30pt];
  \node[dbn boundary,above=3.2pt] at (r1) {$r_1$};
  \node[dbn boundary,above=3.2pt] at (r2) {$r_2$};
  \node[dbn boundary,text=wallred,below=3.2pt] at (q1) {$1$};
  \node[dbn boundary,text=wallred,below=3.2pt] at (q2) {$2$};
  % Local out-port of t1 on arc 3.
  \node[dbn port] (p17-e3-out) at (0.30000,1.40000) {$L$};
  % Local in-port of h1 on arc 3.
  \node[dbn port] (p17-e3-in) at (0.30000,0.70000) {$L$};
  % Local out-port of t1 on arc 4.
  \node[dbn port] (p17-e4-out) at (1.34000,1.40000) {$R$};
  % Local in-port of h1 on arc 4.
  \node[dbn port] (p17-e4-in) at (1.34000,0.70000) {$R$};
\end{tikzpicture}
\dbncaption{17}
\end{minipage}%
\hspace{1.4mm}%
\begin{minipage}[t][3.18cm][t]{0.235\textwidth}\centering
% Class 18: same numbering, vertices, arcs and port orders as the uploaded HTML.
\begin{tikzpicture}[x=0.92cm,y=0.92cm]
  \path[use as bounding box] (-1.60,-0.55) rectangle (1.60,2.60);
  \coordinate (r1) at (-0.82000,1.68000);
  \coordinate (r2) at (0.82000,2.10000);
  \coordinate (t1) at (0.82000,1.50000);
  \coordinate (h1) at (0.82000,0.60000);
  \coordinate (q1) at (-0.82000,0.42000);
  \coordinate (q2) at (0.82000,0.00000);
  % Arc 1: r1[-] -> q1[-].
  \draw[dbn edge] (r1) -- (q1);
  % Arc 2: r2[-] -> t1[-].
  \draw[dbn edge] (r2) -- (t1);
  % Arc 3: t1[L] -> h1[R].
  \draw[dbn edge] (t1) .. controls (0.34000,1.29000) and (0.34000,0.81000) .. (h1);
  % Arc 4: t1[R] -> h1[L].
  \draw[dbn edge] (t1) .. controls (1.30000,1.29000) and (1.30000,0.81000) .. (h1);
  % Arc 5: h1[-] -> q2[-].
  \draw[dbn edge] (h1) -- (q2);
  \fill[black] (r1) circle[radius=2.30pt];
  \fill[black] (r2) circle[radius=2.30pt];
  \fill[black] (t1) circle[radius=2.30pt];
  \fill[black] (h1) circle[radius=2.30pt];
  \fill[black] (q1) circle[radius=2.30pt];
  \fill[black] (q2) circle[radius=2.30pt];
  \node[dbn boundary,above=3.2pt] at (r1) {$r_1$};
  \node[dbn boundary,above=3.2pt] at (r2) {$r_2$};
  \node[dbn boundary,text=wallred,below=3.2pt] at (q1) {$1$};
  \node[dbn boundary,text=wallred,below=3.2pt] at (q2) {$2$};
  % Local out-port of t1 on arc 3.
  \node[dbn port] (p18-e3-out) at (0.30000,1.40000) {$L$};
  % Local in-port of h1 on arc 3.
  \node[dbn port] (p18-e3-in) at (0.30000,0.70000) {$R$};
  % Local out-port of t1 on arc 4.
  \node[dbn port] (p18-e4-out) at (1.34000,1.40000) {$R$};
  % Local in-port of h1 on arc 4.
  \node[dbn port] (p18-e4-in) at (1.34000,0.70000) {$L$};
\end{tikzpicture}
\dbncaption{18}
\end{minipage}%
\hspace{1.4mm}%
\begin{minipage}[t][3.18cm][t]{0.235\textwidth}\centering
% Class 19: same numbering, vertices, arcs and port orders as the uploaded HTML.
\begin{tikzpicture}[x=0.92cm,y=0.92cm]
  \path[use as bounding box] (-1.60,-0.55) rectangle (1.60,2.60);
  \coordinate (r1) at (-1.25000,2.10000);
  \coordinate (r2) at (1.25000,2.10000);
  \coordinate (t1) at (-0.48000,1.38000);
  \coordinate (h1) at (0.48000,0.72000);
  \coordinate (q1) at (-1.25000,0.00000);
  \coordinate (q2) at (1.25000,0.00000);
  % Arc 1: r1[-] -> t1[-].
  \draw[dbn edge] (r1) -- (t1);
  % Arc 2: r2[-] -> h1[L].
  \draw[dbn edge] (r2) -- (h1);
  % Arc 3: t1[R] -> q1[-].
  \draw[dbn edge] (t1) -- (q1);
  % Arc 4: t1[L] -> h1[R].
  \draw[dbn edge] (t1) -- (h1);
  % Arc 5: h1[-] -> q2[-].
  \draw[dbn edge] (h1) -- (q2);
  \fill[black] (r1) circle[radius=2.30pt];
  \fill[black] (r2) circle[radius=2.30pt];
  \fill[black] (t1) circle[radius=2.30pt];
  \fill[black] (h1) circle[radius=2.30pt];
  \fill[black] (q1) circle[radius=2.30pt];
  \fill[black] (q2) circle[radius=2.30pt];
  \node[dbn boundary,above=3.2pt] at (r1) {$r_1$};
  \node[dbn boundary,above=3.2pt] at (r2) {$r_2$};
  \node[dbn boundary,text=wallred,below=3.2pt] at (q1) {$1$};
  \node[dbn boundary,text=wallred,below=3.2pt] at (q2) {$2$};
  % Local in-port of h1 on arc 2.
  \node[dbn port] (p19-e2-in) at (0.88429,0.84941) {$L$};
  % Local out-port of t1 on arc 3.
  \node[dbn port] (p19-e3-out) at (-0.88429,1.25059) {$R$};
  % Local out-port of t1 on arc 4.
  \node[dbn port] (p19-e4-out) at (-0.06025,1.44335) {$L$};
  % Local in-port of h1 on arc 4.
  \node[dbn port] (p19-e4-in) at (0.06025,0.65665) {$R$};
\end{tikzpicture}
\dbncaption{19}
\end{minipage}%
\hspace{1.4mm}%
\begin{minipage}[t][3.18cm][t]{0.235\textwidth}\centering
% Class 20: same numbering, vertices, arcs and port orders as the uploaded HTML.
\begin{tikzpicture}[x=0.92cm,y=0.92cm]
  \path[use as bounding box] (-1.60,-0.55) rectangle (1.60,2.60);
  \coordinate (r1) at (-1.25000,2.10000);
  \coordinate (r2) at (1.25000,2.10000);
  \coordinate (t1) at (-0.48000,1.38000);
  \coordinate (h1) at (0.48000,0.72000);
  \coordinate (q1) at (-1.25000,0.00000);
  \coordinate (q2) at (1.25000,0.00000);
  % Arc 1: r1[-] -> t1[-].
  \draw[dbn edge] (r1) -- (t1);
  % Arc 2: r2[-] -> h1[L].
  \draw[dbn edge] (r2) -- (h1);
  % Arc 3: t1[L] -> q1[-].
  \draw[dbn edge] (t1) -- (q1);
  % Arc 4: t1[R] -> h1[R].
  \draw[dbn edge] (t1) -- (h1);
  % Arc 5: h1[-] -> q2[-].
  \draw[dbn edge] (h1) -- (q2);
  \fill[black] (r1) circle[radius=2.30pt];
  \fill[black] (r2) circle[radius=2.30pt];
  \fill[black] (t1) circle[radius=2.30pt];
  \fill[black] (h1) circle[radius=2.30pt];
  \fill[black] (q1) circle[radius=2.30pt];
  \fill[black] (q2) circle[radius=2.30pt];
  \node[dbn boundary,above=3.2pt] at (r1) {$r_1$};
  \node[dbn boundary,above=3.2pt] at (r2) {$r_2$};
  \node[dbn boundary,text=wallred,below=3.2pt] at (q1) {$1$};
  \node[dbn boundary,text=wallred,below=3.2pt] at (q2) {$2$};
  % Local in-port of h1 on arc 2.
  \node[dbn port] (p20-e2-in) at (0.88429,0.84941) {$L$};
  % Local out-port of t1 on arc 3.
  \node[dbn port] (p20-e3-out) at (-0.88429,1.25059) {$L$};
  % Local out-port of t1 on arc 4.
  \node[dbn port] (p20-e4-out) at (-0.06025,1.44335) {$R$};
  % Local in-port of h1 on arc 4.
  \node[dbn port] (p20-e4-in) at (0.06025,0.65665) {$R$};
\end{tikzpicture}
\dbncaption{20}
\end{minipage}%
}\par

\noindent\makebox[\textwidth][c]{%
\begin{minipage}[t][3.18cm][t]{0.235\textwidth}\centering
% Class 21: same numbering, vertices, arcs and port orders as the uploaded HTML.
\begin{tikzpicture}[x=0.92cm,y=0.92cm]
  \path[use as bounding box] (-1.60,-0.55) rectangle (1.60,2.60);
  \coordinate (r1) at (-1.25000,2.10000);
  \coordinate (r2) at (1.25000,2.10000);
  \coordinate (t1) at (-0.48000,1.38000);
  \coordinate (h1) at (0.48000,0.72000);
  \coordinate (q1) at (-1.25000,0.00000);
  \coordinate (q2) at (1.25000,0.00000);
  % Arc 1: r1[-] -> t1[-].
  \draw[dbn edge] (r1) -- (t1);
  % Arc 2: r2[-] -> h1[R].
  \draw[dbn edge] (r2) -- (h1);
  % Arc 3: t1[R] -> q1[-].
  \draw[dbn edge] (t1) -- (q1);
  % Arc 4: t1[L] -> h1[L].
  \draw[dbn edge] (t1) -- (h1);
  % Arc 5: h1[-] -> q2[-].
  \draw[dbn edge] (h1) -- (q2);
  \fill[black] (r1) circle[radius=2.30pt];
  \fill[black] (r2) circle[radius=2.30pt];
  \fill[black] (t1) circle[radius=2.30pt];
  \fill[black] (h1) circle[radius=2.30pt];
  \fill[black] (q1) circle[radius=2.30pt];
  \fill[black] (q2) circle[radius=2.30pt];
  \node[dbn boundary,above=3.2pt] at (r1) {$r_1$};
  \node[dbn boundary,above=3.2pt] at (r2) {$r_2$};
  \node[dbn boundary,text=wallred,below=3.2pt] at (q1) {$1$};
  \node[dbn boundary,text=wallred,below=3.2pt] at (q2) {$2$};
  % Local in-port of h1 on arc 2.
  \node[dbn port] (p21-e2-in) at (0.88429,0.84941) {$R$};
  % Local out-port of t1 on arc 3.
  \node[dbn port] (p21-e3-out) at (-0.88429,1.25059) {$R$};
  % Local out-port of t1 on arc 4.
  \node[dbn port] (p21-e4-out) at (-0.06025,1.44335) {$L$};
  % Local in-port of h1 on arc 4.
  \node[dbn port] (p21-e4-in) at (0.06025,0.65665) {$L$};
\end{tikzpicture}
\dbncaption{21}
\end{minipage}%
\hspace{1.4mm}%
\begin{minipage}[t][3.18cm][t]{0.235\textwidth}\centering
% Class 22: same numbering, vertices, arcs and port orders as the uploaded HTML.
\begin{tikzpicture}[x=0.92cm,y=0.92cm]
  \path[use as bounding box] (-1.60,-0.55) rectangle (1.60,2.60);
  \coordinate (r1) at (-1.25000,2.10000);
  \coordinate (r2) at (1.25000,2.10000);
  \coordinate (t1) at (-0.48000,1.38000);
  \coordinate (h1) at (0.48000,0.72000);
  \coordinate (q1) at (-1.25000,0.00000);
  \coordinate (q2) at (1.25000,0.00000);
  % Arc 1: r1[-] -> t1[-].
  \draw[dbn edge] (r1) -- (t1);
  % Arc 2: r2[-] -> h1[R].
  \draw[dbn edge] (r2) -- (h1);
  % Arc 3: t1[L] -> q1[-].
  \draw[dbn edge] (t1) -- (q1);
  % Arc 4: t1[R] -> h1[L].
  \draw[dbn edge] (t1) -- (h1);
  % Arc 5: h1[-] -> q2[-].
  \draw[dbn edge] (h1) -- (q2);
  \fill[black] (r1) circle[radius=2.30pt];
  \fill[black] (r2) circle[radius=2.30pt];
  \fill[black] (t1) circle[radius=2.30pt];
  \fill[black] (h1) circle[radius=2.30pt];
  \fill[black] (q1) circle[radius=2.30pt];
  \fill[black] (q2) circle[radius=2.30pt];
  \node[dbn boundary,above=3.2pt] at (r1) {$r_1$};
  \node[dbn boundary,above=3.2pt] at (r2) {$r_2$};
  \node[dbn boundary,text=wallred,below=3.2pt] at (q1) {$1$};
  \node[dbn boundary,text=wallred,below=3.2pt] at (q2) {$2$};
  % Local in-port of h1 on arc 2.
  \node[dbn port] (p22-e2-in) at (0.88429,0.84941) {$R$};
  % Local out-port of t1 on arc 3.
  \node[dbn port] (p22-e3-out) at (-0.88429,1.25059) {$L$};
  % Local out-port of t1 on arc 4.
  \node[dbn port] (p22-e4-out) at (-0.06025,1.44335) {$R$};
  % Local in-port of h1 on arc 4.
  \node[dbn port] (p22-e4-in) at (0.06025,0.65665) {$L$};
\end{tikzpicture}
\dbncaption{22}
\end{minipage}%
\hspace{1.4mm}%
\begin{minipage}[t][3.18cm][t]{0.235\textwidth}\centering
% Class 23: same numbering, vertices, arcs and port orders as the uploaded HTML.
\begin{tikzpicture}[x=0.92cm,y=0.92cm]
  \path[use as bounding box] (-1.60,-0.55) rectangle (1.60,2.60);
  \coordinate (r1) at (-0.82000,2.10000);
  \coordinate (r2) at (0.82000,1.68000);
  \coordinate (t1) at (-0.82000,1.50000);
  \coordinate (h1) at (-0.82000,0.60000);
  \coordinate (q1) at (0.82000,0.42000);
  \coordinate (q2) at (-0.82000,0.00000);
  % Arc 1: r1[-] -> t1[-].
  \draw[dbn edge] (r1) -- (t1);
  % Arc 2: r2[-] -> q1[-].
  \draw[dbn edge] (r2) -- (q1);
  % Arc 3: t1[L] -> h1[L].
  \draw[dbn edge] (t1) .. controls (-1.30000,1.29000) and (-1.30000,0.81000) .. (h1);
  % Arc 4: t1[R] -> h1[R].
  \draw[dbn edge] (t1) .. controls (-0.34000,1.29000) and (-0.34000,0.81000) .. (h1);
  % Arc 5: h1[-] -> q2[-].
  \draw[dbn edge] (h1) -- (q2);
  \fill[black] (r1) circle[radius=2.30pt];
  \fill[black] (r2) circle[radius=2.30pt];
  \fill[black] (t1) circle[radius=2.30pt];
  \fill[black] (h1) circle[radius=2.30pt];
  \fill[black] (q1) circle[radius=2.30pt];
  \fill[black] (q2) circle[radius=2.30pt];
  \node[dbn boundary,above=3.2pt] at (r1) {$r_1$};
  \node[dbn boundary,above=3.2pt] at (r2) {$r_2$};
  \node[dbn boundary,text=wallred,below=3.2pt] at (q1) {$1$};
  \node[dbn boundary,text=wallred,below=3.2pt] at (q2) {$2$};
  % Local out-port of t1 on arc 3.
  \node[dbn port] (p23-e3-out) at (-1.34000,1.40000) {$L$};
  % Local in-port of h1 on arc 3.
  \node[dbn port] (p23-e3-in) at (-1.34000,0.70000) {$L$};
  % Local out-port of t1 on arc 4.
  \node[dbn port] (p23-e4-out) at (-0.30000,1.40000) {$R$};
  % Local in-port of h1 on arc 4.
  \node[dbn port] (p23-e4-in) at (-0.30000,0.70000) {$R$};
\end{tikzpicture}
\dbncaption{23}
\end{minipage}%
\hspace{1.4mm}%
\begin{minipage}[t][3.18cm][t]{0.235\textwidth}\centering
% Class 24: same numbering, vertices, arcs and port orders as the uploaded HTML.
\begin{tikzpicture}[x=0.92cm,y=0.92cm]
  \path[use as bounding box] (-1.60,-0.55) rectangle (1.60,2.60);
  \coordinate (r1) at (-0.82000,2.10000);
  \coordinate (r2) at (0.82000,1.68000);
  \coordinate (t1) at (-0.82000,1.50000);
  \coordinate (h1) at (-0.82000,0.60000);
  \coordinate (q1) at (0.82000,0.42000);
  \coordinate (q2) at (-0.82000,0.00000);
  % Arc 1: r1[-] -> t1[-].
  \draw[dbn edge] (r1) -- (t1);
  % Arc 2: r2[-] -> q1[-].
  \draw[dbn edge] (r2) -- (q1);
  % Arc 3: t1[L] -> h1[R].
  \draw[dbn edge] (t1) .. controls (-1.30000,1.29000) and (-1.30000,0.81000) .. (h1);
  % Arc 4: t1[R] -> h1[L].
  \draw[dbn edge] (t1) .. controls (-0.34000,1.29000) and (-0.34000,0.81000) .. (h1);
  % Arc 5: h1[-] -> q2[-].
  \draw[dbn edge] (h1) -- (q2);
  \fill[black] (r1) circle[radius=2.30pt];
  \fill[black] (r2) circle[radius=2.30pt];
  \fill[black] (t1) circle[radius=2.30pt];
  \fill[black] (h1) circle[radius=2.30pt];
  \fill[black] (q1) circle[radius=2.30pt];
  \fill[black] (q2) circle[radius=2.30pt];
  \node[dbn boundary,above=3.2pt] at (r1) {$r_1$};
  \node[dbn boundary,above=3.2pt] at (r2) {$r_2$};
  \node[dbn boundary,text=wallred,below=3.2pt] at (q1) {$1$};
  \node[dbn boundary,text=wallred,below=3.2pt] at (q2) {$2$};
  % Local out-port of t1 on arc 3.
  \node[dbn port] (p24-e3-out) at (-1.34000,1.40000) {$L$};
  % Local in-port of h1 on arc 3.
  \node[dbn port] (p24-e3-in) at (-1.34000,0.70000) {$R$};
  % Local out-port of t1 on arc 4.
  \node[dbn port] (p24-e4-out) at (-0.30000,1.40000) {$R$};
  % Local in-port of h1 on arc 4.
  \node[dbn port] (p24-e4-in) at (-0.30000,0.70000) {$L$};
\end{tikzpicture}
\dbncaption{24}
\end{minipage}%
}\par

\noindent\makebox[\textwidth][c]{%
\begin{minipage}[t][3.18cm][t]{0.235\textwidth}\centering
% Class 25: same numbering, vertices, arcs and port orders as the uploaded HTML.
\begin{tikzpicture}[x=0.92cm,y=0.92cm]
  \path[use as bounding box] (-1.60,-0.55) rectangle (1.60,2.60);
  \coordinate (r1) at (-1.05000,2.10000);
  \coordinate (r2) at (1.05000,2.10000);
  \coordinate (t1) at (0.00000,0.70000);
  \coordinate (h1) at (0.00000,1.40000);
  \coordinate (q1) at (-1.05000,0.00000);
  \coordinate (q2) at (1.05000,0.00000);
  % Arc 1: r1[-] -> h1[L].
  \draw[dbn edge] (r1) -- (h1);
  % Arc 2: r2[-] -> h1[R].
  \draw[dbn edge] (r2) -- (h1);
  % Arc 3: t1[L] -> q1[-].
  \draw[dbn edge] (t1) -- (q1);
  % Arc 4: t1[R] -> q2[-].
  \draw[dbn edge] (t1) -- (q2);
  % Arc 5: h1[-] -> t1[-].
  \draw[dbn edge] (h1) -- (t1);
  \fill[black] (r1) circle[radius=2.30pt];
  \fill[black] (r2) circle[radius=2.30pt];
  \fill[black] (t1) circle[radius=2.30pt];
  \fill[black] (h1) circle[radius=2.30pt];
  \fill[black] (q1) circle[radius=2.30pt];
  \fill[black] (q2) circle[radius=2.30pt];
  \node[dbn boundary,above=3.2pt] at (r1) {$r_1$};
  \node[dbn boundary,above=3.2pt] at (r2) {$r_2$};
  \node[dbn boundary,text=wallred,below=3.2pt] at (q1) {$1$};
  \node[dbn boundary,text=wallred,below=3.2pt] at (q2) {$2$};
  % Local in-port of h1 on arc 1.
  \node[dbn port] (p25-e1-in) at (-0.20000,1.87000) {$L$};
  % Local in-port of h1 on arc 2.
  \node[dbn port] (p25-e2-in) at (0.20000,1.87000) {$R$};
  % Local out-port of t1 on arc 3.
  \node[dbn port] (p25-e3-out) at (-0.20000,0.23000) {$L$};
  % Local out-port of t1 on arc 4.
  \node[dbn port] (p25-e4-out) at (0.20000,0.23000) {$R$};
\end{tikzpicture}
\dbncaption{25}
\end{minipage}%
\hspace{1.4mm}%
\begin{minipage}[t][3.18cm][t]{0.235\textwidth}\centering
% Class 26: same numbering, vertices, arcs and port orders as the uploaded HTML.
\begin{tikzpicture}[x=0.92cm,y=0.92cm]
  \path[use as bounding box] (-1.60,-0.55) rectangle (1.60,2.60);
  \coordinate (r1) at (-1.05000,2.10000);
  \coordinate (r2) at (1.05000,2.10000);
  \coordinate (t1) at (0.00000,0.70000);
  \coordinate (h1) at (0.00000,1.40000);
  \coordinate (q1) at (-1.05000,0.00000);
  \coordinate (q2) at (1.05000,0.00000);
  % Arc 1: r1[-] -> h1[L].
  \draw[dbn edge] (r1) -- (h1);
  % Arc 2: r2[-] -> h1[R].
  \draw[dbn edge] (r2) -- (h1);
  % Arc 3: t1[R] -> q1[-].
  \draw[dbn edge] (t1) -- (q1);
  % Arc 4: t1[L] -> q2[-].
  \draw[dbn edge] (t1) -- (q2);
  % Arc 5: h1[-] -> t1[-].
  \draw[dbn edge] (h1) -- (t1);
  \fill[black] (r1) circle[radius=2.30pt];
  \fill[black] (r2) circle[radius=2.30pt];
  \fill[black] (t1) circle[radius=2.30pt];
  \fill[black] (h1) circle[radius=2.30pt];
  \fill[black] (q1) circle[radius=2.30pt];
  \fill[black] (q2) circle[radius=2.30pt];
  \node[dbn boundary,above=3.2pt] at (r1) {$r_1$};
  \node[dbn boundary,above=3.2pt] at (r2) {$r_2$};
  \node[dbn boundary,text=wallred,below=3.2pt] at (q1) {$1$};
  \node[dbn boundary,text=wallred,below=3.2pt] at (q2) {$2$};
  % Local in-port of h1 on arc 1.
  \node[dbn port] (p26-e1-in) at (-0.20000,1.87000) {$L$};
  % Local in-port of h1 on arc 2.
  \node[dbn port] (p26-e2-in) at (0.20000,1.87000) {$R$};
  % Local out-port of t1 on arc 3.
  \node[dbn port] (p26-e3-out) at (-0.20000,0.23000) {$R$};
  % Local out-port of t1 on arc 4.
  \node[dbn port] (p26-e4-out) at (0.20000,0.23000) {$L$};
\end{tikzpicture}
\dbncaption{26}
\end{minipage}%
\hspace{1.4mm}%
\begin{minipage}[t][3.18cm][t]{0.235\textwidth}\centering
% Class 27: same numbering, vertices, arcs and port orders as the uploaded HTML.
\begin{tikzpicture}[x=0.92cm,y=0.92cm]
  \path[use as bounding box] (-1.60,-0.55) rectangle (1.60,2.60);
  \coordinate (r1) at (-1.05000,2.10000);
  \coordinate (r2) at (1.05000,2.10000);
  \coordinate (t1) at (0.00000,0.70000);
  \coordinate (h1) at (0.00000,1.40000);
  \coordinate (q1) at (-1.05000,0.00000);
  \coordinate (q2) at (1.05000,0.00000);
  % Arc 1: r1[-] -> h1[R].
  \draw[dbn edge] (r1) -- (h1);
  % Arc 2: r2[-] -> h1[L].
  \draw[dbn edge] (r2) -- (h1);
  % Arc 3: t1[L] -> q1[-].
  \draw[dbn edge] (t1) -- (q1);
  % Arc 4: t1[R] -> q2[-].
  \draw[dbn edge] (t1) -- (q2);
  % Arc 5: h1[-] -> t1[-].
  \draw[dbn edge] (h1) -- (t1);
  \fill[black] (r1) circle[radius=2.30pt];
  \fill[black] (r2) circle[radius=2.30pt];
  \fill[black] (t1) circle[radius=2.30pt];
  \fill[black] (h1) circle[radius=2.30pt];
  \fill[black] (q1) circle[radius=2.30pt];
  \fill[black] (q2) circle[radius=2.30pt];
  \node[dbn boundary,above=3.2pt] at (r1) {$r_1$};
  \node[dbn boundary,above=3.2pt] at (r2) {$r_2$};
  \node[dbn boundary,text=wallred,below=3.2pt] at (q1) {$1$};
  \node[dbn boundary,text=wallred,below=3.2pt] at (q2) {$2$};
  % Local in-port of h1 on arc 1.
  \node[dbn port] (p27-e1-in) at (-0.20000,1.87000) {$R$};
  % Local in-port of h1 on arc 2.
  \node[dbn port] (p27-e2-in) at (0.20000,1.87000) {$L$};
  % Local out-port of t1 on arc 3.
  \node[dbn port] (p27-e3-out) at (-0.20000,0.23000) {$L$};
  % Local out-port of t1 on arc 4.
  \node[dbn port] (p27-e4-out) at (0.20000,0.23000) {$R$};
\end{tikzpicture}
\dbncaption{27}
\end{minipage}%
\hspace{1.4mm}%
\begin{minipage}[t][3.18cm][t]{0.235\textwidth}\centering
% Class 28: same numbering, vertices, arcs and port orders as the uploaded HTML.
\begin{tikzpicture}[x=0.92cm,y=0.92cm]
  \path[use as bounding box] (-1.60,-0.55) rectangle (1.60,2.60);
  \coordinate (r1) at (-1.05000,2.10000);
  \coordinate (r2) at (1.05000,2.10000);
  \coordinate (t1) at (0.00000,0.70000);
  \coordinate (h1) at (0.00000,1.40000);
  \coordinate (q1) at (-1.05000,0.00000);
  \coordinate (q2) at (1.05000,0.00000);
  % Arc 1: r1[-] -> h1[R].
  \draw[dbn edge] (r1) -- (h1);
  % Arc 2: r2[-] -> h1[L].
  \draw[dbn edge] (r2) -- (h1);
  % Arc 3: t1[R] -> q1[-].
  \draw[dbn edge] (t1) -- (q1);
  % Arc 4: t1[L] -> q2[-].
  \draw[dbn edge] (t1) -- (q2);
  % Arc 5: h1[-] -> t1[-].
  \draw[dbn edge] (h1) -- (t1);
  \fill[black] (r1) circle[radius=2.30pt];
  \fill[black] (r2) circle[radius=2.30pt];
  \fill[black] (t1) circle[radius=2.30pt];
  \fill[black] (h1) circle[radius=2.30pt];
  \fill[black] (q1) circle[radius=2.30pt];
  \fill[black] (q2) circle[radius=2.30pt];
  \node[dbn boundary,above=3.2pt] at (r1) {$r_1$};
  \node[dbn boundary,above=3.2pt] at (r2) {$r_2$};
  \node[dbn boundary,text=wallred,below=3.2pt] at (q1) {$1$};
  \node[dbn boundary,text=wallred,below=3.2pt] at (q2) {$2$};
  % Local in-port of h1 on arc 1.
  \node[dbn port] (p28-e1-in) at (-0.20000,1.87000) {$R$};
  % Local in-port of h1 on arc 2.
  \node[dbn port] (p28-e2-in) at (0.20000,1.87000) {$L$};
  % Local out-port of t1 on arc 3.
  \node[dbn port] (p28-e3-out) at (-0.20000,0.23000) {$R$};
  % Local out-port of t1 on arc 4.
  \node[dbn port] (p28-e4-out) at (0.20000,0.23000) {$L$};
\end{tikzpicture}
\dbncaption{28}
\end{minipage}%
}\par

\newcommand{\DecoratedTwo}[8]{%
\begin{tikzpicture}[
  x=0.80cm,y=0.70cm,
  line width=0.43pt,
  every node/.style={font=\small,inner sep=0pt,outer sep=0pt},
  dt hole/.style={draw=holegray,fill=holefill,dash pattern=on 2pt off 1.4pt,line width=0.43pt},
  dt wall/.style={draw=wallred,line width=0.95pt}
]
  \path[use as bounding box] (-0.12,-0.13) rectangle (3.06,3.09);
  % Top and middle rows.
  \draw (0,2) rectangle (1,3);
  \draw (0,1) rectangle (1,2);
  \draw (1,1) rectangle (2,2);
  \node at (0.5,2.5) {$#1$};
  \node at (0.5,1.5) {$#2$};
  \node at (1.5,1.5) {$#3$};
  % Decorated bottom row.
  \ifnum#4=1\relax
    \draw (0,0) rectangle (1,1);
    \draw[dt hole] (1,0) rectangle (2,1);
    \node at (0.50,0.55) {$#5$};
    \node[font=\fontsize{7}{8}\selectfont,anchor=south east]
       at (0.93,0.08) {$#6$};
    \node[text=wallred] at (1.5,0.5) {$#7$};
  \else
    \draw[dt hole] (0,0) rectangle (1,1);
    \draw (1,0) rectangle (2,1);
    \node[text=wallred] at (0.5,0.5) {$#7$};
    \node at (1.50,0.55) {$#5$};
    \node[font=\fontsize{7}{8}\selectfont,anchor=south east]
       at (1.93,0.08) {$#6$};
  \fi
  \draw[dt wall] (1,0.06)--(1,0.94);
  % Star decoration.
  \node at (2.27,0.50) {$\star$};
  \node[text=wallred] at (2.75,0.50) {$#8$};
\end{tikzpicture}%
}

\section{Decorated Young tableaux of type $(2,1,1)$}
\label{sec:AppendixB}
\begin{center}
  {\small Corresponding to decorated binary networks (1)--(28) in Appendix A under $\Psi$}\par
\end{center}
\vspace{1mm}
\setlength{\parindent}{0pt}
\newcommand{\TableauPanel}[2]{%
\begin{minipage}[t][3.20cm][t]{0.246\textwidth}
  \centering
  #2\par
  \vspace{0.5mm}
  {\footnotesize (#1)}\par
\end{minipage}%
}
\noindent\makebox[\textwidth][c]{%
\TableauPanel{1}{\DecoratedTwo{4}{2}{3}{1}{1}{R}{2}{1}}%
\hspace{0.4mm}%
\TableauPanel{2}{\DecoratedTwo{4}{2}{3}{2}{1}{R}{2}{1}}%
\hspace{0.4mm}%
\TableauPanel{3}{\DecoratedTwo{4}{2}{3}{1}{1}{L}{2}{1}}%
\hspace{0.4mm}%
\TableauPanel{4}{\DecoratedTwo{4}{2}{3}{2}{1}{L}{2}{1}}%
}\par
\noindent\makebox[\textwidth][c]{%
\TableauPanel{5}{\DecoratedTwo{4}{1}{3}{2}{2}{R}{1}{2}}%
\hspace{0.4mm}%
\TableauPanel{6}{\DecoratedTwo{4}{1}{3}{2}{2}{L}{1}{2}}%
\hspace{0.4mm}%
\TableauPanel{7}{\DecoratedTwo{3}{2}{4}{2}{1}{L}{2}{1}}%
\hspace{0.4mm}%
\TableauPanel{8}{\DecoratedTwo{3}{1}{4}{2}{2}{L}{1}{2}}%
}\par
\noindent\makebox[\textwidth][c]{%
\TableauPanel{9}{\DecoratedTwo{3}{2}{4}{2}{1}{R}{2}{1}}%
\hspace{0.4mm}%
\TableauPanel{10}{\DecoratedTwo{3}{1}{4}{2}{2}{R}{1}{2}}%
\hspace{0.4mm}%
\TableauPanel{11}{\DecoratedTwo{3}{2}{4}{1}{1}{R}{2}{1}}%
\hspace{0.4mm}%
\TableauPanel{12}{\DecoratedTwo{3}{2}{4}{1}{1}{L}{2}{1}}%
}\par
\noindent\makebox[\textwidth][c]{%
\TableauPanel{13}{\DecoratedTwo{4}{2}{3}{1}{1}{R}{1}{2}}%
\hspace{0.4mm}%
\TableauPanel{14}{\DecoratedTwo{4}{2}{3}{2}{1}{R}{1}{2}}%
\hspace{0.4mm}%
\TableauPanel{15}{\DecoratedTwo{4}{2}{3}{1}{1}{L}{1}{2}}%
\hspace{0.4mm}%
\TableauPanel{16}{\DecoratedTwo{4}{2}{3}{2}{1}{L}{1}{2}}%
}\par
\noindent\makebox[\textwidth][c]{%
\TableauPanel{17}{\DecoratedTwo{4}{1}{3}{2}{2}{R}{2}{1}}%
\hspace{0.4mm}%
\TableauPanel{18}{\DecoratedTwo{4}{1}{3}{2}{2}{L}{2}{1}}%
\hspace{0.4mm}%
\TableauPanel{19}{\DecoratedTwo{3}{2}{4}{2}{1}{L}{1}{2}}%
\hspace{0.4mm}%
\TableauPanel{20}{\DecoratedTwo{3}{1}{4}{2}{2}{L}{2}{1}}%
}\par
\noindent\makebox[\textwidth][c]{%
\TableauPanel{21}{\DecoratedTwo{3}{2}{4}{2}{1}{R}{1}{2}}%
\hspace{0.4mm}%
\TableauPanel{22}{\DecoratedTwo{3}{1}{4}{2}{2}{R}{2}{1}}%
\hspace{0.4mm}%
\TableauPanel{23}{\DecoratedTwo{3}{2}{4}{1}{1}{R}{1}{2}}%
\hspace{0.4mm}%
\TableauPanel{24}{\DecoratedTwo{3}{2}{4}{1}{1}{L}{1}{2}}%
}\par
\noindent\makebox[\textwidth][c]{%
\TableauPanel{25}{\DecoratedTwo{2}{1}{4}{2}{3}{R}{2}{1}}%
\hspace{0.4mm}%
\TableauPanel{26}{\DecoratedTwo{2}{1}{4}{2}{3}{R}{1}{2}}%
\hspace{0.4mm}%
\TableauPanel{27}{\DecoratedTwo{2}{1}{4}{2}{3}{L}{2}{1}}%
\hspace{0.4mm}%
\TableauPanel{28}{\DecoratedTwo{2}{1}{4}{2}{3}{L}{1}{2}}%
}\par

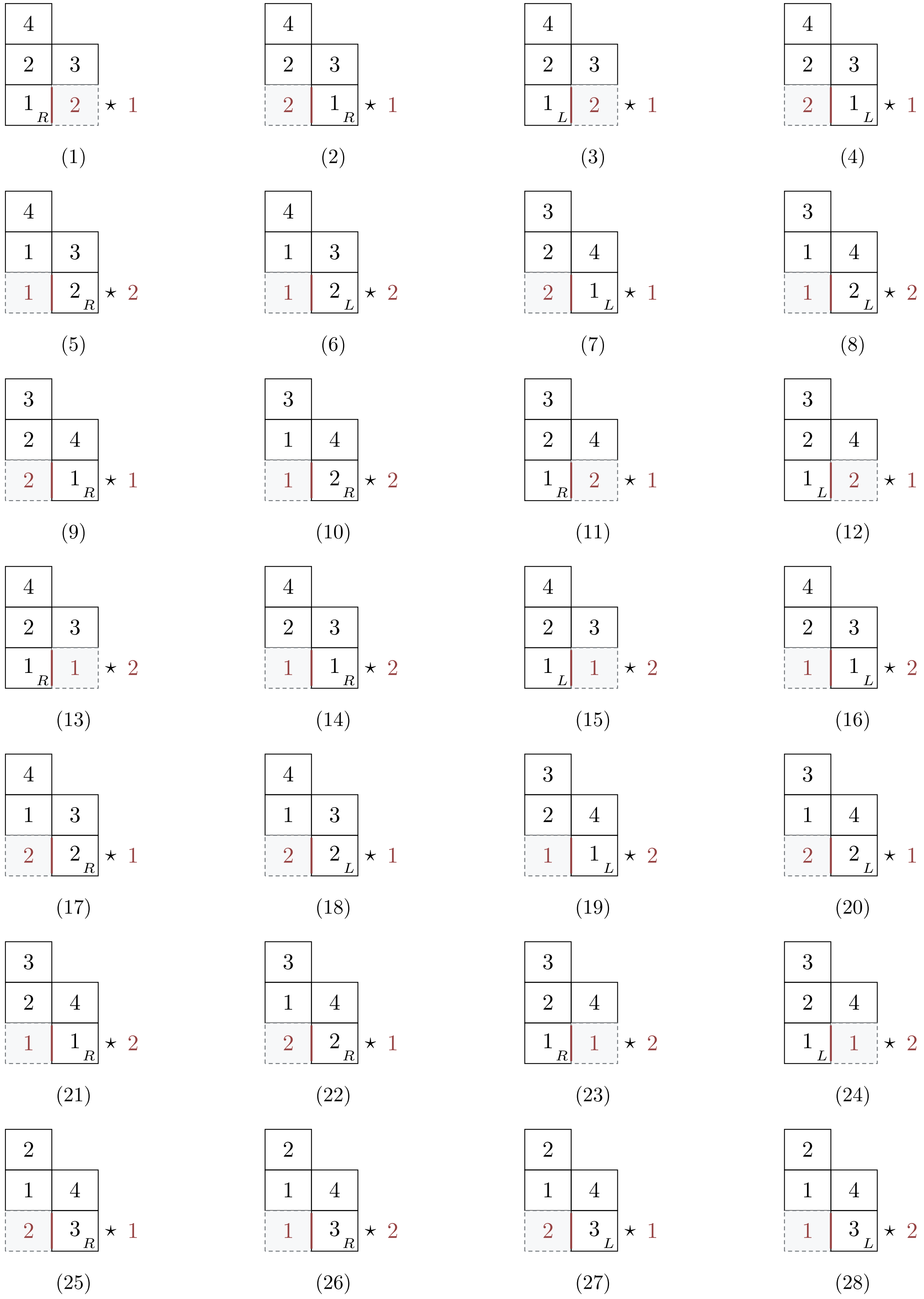
\end{document}